\documentclass[12pt,leqno]{amsart}  
\usepackage{amsmath,amstext,amsthm,amssymb,amsxtra}
\usepackage[colorlinks,citecolor=red,pagebackref,hypertexnames=false]{hyperref}
\usepackage[backrefs]{amsrefs}
\allowdisplaybreaks
\swapnumbers
 
\theoremstyle{plain} 
\newtheorem{lemma}[equation]{Lemma} 
\newtheorem{proposition}[equation]{Proposition} 
\newtheorem{theorem}[equation]{Theorem} 
\newtheorem{corollary}[equation]{Corollary} 
\newtheorem*{nehari}{Nehari Theorem}
\newtheorem*{multinehari}{Multiparameter Nehari Theorem}

\theoremstyle{definition}

\theoremstyle{remark}
\newtheorem{remark}[equation]{Remark}
\newtheorem*{acknowledgment}{Acknowledgment}

\numberwithin{equation}{section}

\def\norm#1.#2.{\lVert#1\rVert_{#2}}
\def\Norm#1.#2.{\bigl\lVert#1\bigr\rVert_{#2}}
\def\NOrm#1.#2.{\Bigl\lVert#1\Bigr\rVert_{#2}}
\def\NORm#1.#2.{\biggl\lVert#1\biggr\rVert_{#2}}
\def\NORM#1.#2.{\Biggl\lVert#1\Biggr\rVert_{#2}}

\def\ip#1,#2,{\langle #1,#2\rangle}
\def\Ip#1,#2,{\bigl\langle#1,#2\bigr\rangle}
\def\IP#1,#2,{\Bigl\langle#1,#2\Bigr\rangle}

\def\mid{\,:\,}

\def\abs#1{\lvert#1\rvert}
\def\Abs#1{\bigl\lvert#1\bigr\rvert}

\def\XXint#1#2#3{{\setbox0=\hbox{$#1{#2#3}{\int}$}
     \vcenter{\hbox{$#2#3$}}\kern-.5\wd0}}

\def\eqdef{\stackrel{\mathrm{def}}{{}={}}}

\begin{document}
\title {Lectures on Nehari's Theorem on the Polydisk}

\author[M.T. Lacey]{Michael T. Lacey\\
School of Mathematics\\
Georgia Institute of Technology\\
Atlanta GA 30332\\
}

\thanks{Research supported in part by a National 
Science Foundation Grant. The author is a Guggenheim Fellow.}

\email{lacey@math.gatech.edu}

\begin{abstract}
We are concerned with Nehari's theorem on Hardy space on a polydisk. 
Define `little' Hankel operators on product  Hardy space $ H^2 ( \mathbb C _+ ^{d})$
by 
\begin{equation*}
\operatorname H _{b} \varphi  \eqdef \operatorname P _{\oplus}  \operatorname M _{b} 
\overline \varphi \,. 
\end{equation*}
where $ \operatorname P _{\oplus}$ is the orthogonal projection from $ L^2 (\mathbb R ^{d})$
to $ H^2 ( \mathbb C _+ ^{d})$ and $ \operatorname M _{b}$ is the operator of 
multiplication by $ b$.  We present the proof of Ferguson and Lacey \cite {MR1961195} 
and Lacey and Terwelleger \cite{witherin} that we have the equivalence of norms 
\begin{equation*}
\norm H_b. .\cdot \simeq \norm b. BMO (\mathbb C _+ ^{d}).
\end{equation*}
for analytic functions $ b$.  Here, $ BMO (\mathbb C _+ ^{d})$ is 
the dual to $ H^1 (\mathbb C _+ ^{d})$ as discovered by Chang and R.~Fefferman.  
This article begins with the classical Nehari theorem, and presents 
the  necessary background for the proof of the extension above.  
The proof of the extension 
is an induction on parameters, with a bootstrapping argument.  Some of the 
more technical details of the earlier proofs are now seen as consequences of a 
paraproduct theory.   
\end{abstract}

\maketitle

\tableofcontents

 \section{Introduction}
 
These notes concern the subject of Nehari's theorem, on Hardy space of the disk, 
and products 
of the disk.  The theorem on the disk is classical, with three different approaches 
possible; the same question on products of the disk, the polydisk of the title,
is a new result of the author, Sarah 
Ferguson and Erin Terwelleger \cites{MR1961195,witherin}.  
The proof in the product setting is much more complicated, 
with currently only one proof known.  It relies upon a delicate
bootstrapping argument with an induction 
on parameters.  These elements are suggested by the 
harmonic analysis associated with  product theory, as developed by 
S.-Y.~Chang, R.~Fefferman and J.-L.~Journ\'e \cites{MR86g:42038,MR82a:32009,chang,
cf1,
cf2,
MR87g:42028,
MR88d:42028}.  These notes will provide an 
approach to this result that is more leisurely than the research articles on the subject.
We in particular include a great many references, and a description of 
related results and concepts.  The proof of the main theorem we give 
 is a little more `structural' in that the main technical estimates 
are seen as consequences of a theory of paraproducts. 

The key concepts of this paper concern the intertwined topics of Hankel operators, 
Hardy space, Hilbert transforms, commutators, and paraproducts.  Let us describe 
Hankel matrices.  

Consider a function $ b$ on $ L^2(\mathbb T )$, and the operator $ \operatorname M _{b} $ 
of pointwise multiplication by $b$.  That is, $ \operatorname M _{b} \varphi \eqdef b\cdot \varphi $. 
Now, $ L^2 (\mathbb T ) $ has the exponential basis.  We view the circle as embedded in the 
natural way in the complex plane, so that a relevant  basis is $ \{z ^{n}\mid n\in \mathbb Z \}$. 
The decomposition of functions in this basis of course generates the Fourier transform: 
\begin{equation*}
\widehat f(n) \eqdef \int _{\mathbb T } f(z) z ^{-n}\; \abs{ dz}  
\end{equation*}
In the exponential basis, $ \operatorname M _{b}$ has a matrix form, 
\begin{equation*}
\operatorname M_b \longleftrightarrow  \{ \widehat b(i-j)\mid i,j\in \mathbb Z \}.  
\end{equation*}
Restrictions of this matrix give \emph{Hankel} and \emph{Toeplitz} operators.  

The restriction of the matrix to the upper quadrant $ \mathbb N \times \mathbb N$ 
is a \emph{Toeplitz} matrix.  Namely, $ \operatorname T=\{t _{ij}\mid i,j\in \mathbb N\}$ 
is a Toeplitz matrix iff $ t _{ij}=\alpha _{i-j}$ for some  
numerical sequence $ \alpha $ on $ \mathbb Z $. 
In this note, we are principally interested in the boundedness properties of 
operators.  It is easy to see that Toeplitz matrix is bounded iff 
the sequence $ \alpha _j$ are the Fourier coefficients of a bounded function.  
How this statement changes for \emph{Hankel} matrices occupies our attention.

Restricting the matrix to say the quadrant $ \mathbb N\times (-\mathbb N)$ gives a \emph{Hankel} 
matrix.  Namely, $ \operatorname H=\{h _{ij}\mid i,j\in \mathbb N\}$ is a \emph{Hankel} 
matrix iff $ h _{ij}=\alpha _{i+j}$ for some numerical sequence on $ \mathbb N$.

In passing to these restrictions of the matrix for $ \operatorname M _{b}$, we are 
implicitly restricting the matrix on  $ \ell ^2(\mathbb Z ) $ to  one on 
$ \ell ^2 (\mathbb N)$.  Namely, 
a Hankel and Toeplitz matrix are operators on $ \ell^2(\mathbb N)$.   
The natural analog in $ L^2(\mathbb T )$ is \emph{Hardy space} $ H^2(\mathbb T )$. 
By definition, $ H^2(\mathbb T )=H^2 _{+}(\mathbb T )$ is the closed subspace of $ L^2(\mathbb T )$  
generated by $ \{z^n\mid n\ge0\}$.   It is natural to call these functions 
\emph{analytic}, as $ f\in H^2(\mathbb T )$ admit an analytic extension to the 
disk $ \mathbb D $ given by 
\begin{equation*}
F(z)=\sum_{n\ge0} \widehat f(n) z ^{n}\,. 
\end{equation*}
Functions in $ H ^2 _{-}(\mathbb T )=L^2(\mathbb T )\ominus H^2(\mathbb T )$ 
are referred to as \emph {antianalytic}. 

Let us describe the Hankel operators on $ H^2(\mathbb T )$.  Let $ \operatorname P_ \pm$
be the orthogonal projection from $ L^2(\mathbb T )$ onto the subspace $ H^2 _{\pm}(\mathbb T )$. 
A \emph{Hankel operator with symbol $ b$} is an operator 
$ \operatorname H_b$ from $ H^2_+(\mathbb T )$ to 
$ H^2 _{-}(\mathbb T )$ given by $ \operatorname H_b \varphi 
\eqdef \operatorname  P_- \operatorname M_b \overline {\varphi }$. 
It is clear that this definition only depends on the analytic part of $ b$. 

\begin{remark}\label{r.conjugate} The placement of the conjugate symbol 
is somewhat arbitrary, and is adopted in this way only for convenience. 
Richard Rochberg avoids such complications by defining Hankel operators as 
bilinear operators $ \operatorname B$ from $ H^2\times H^2$ into $ \mathbb C$, 
which are linear on products:  $ \operatorname B(\varphi ,\psi )=\operatorname L(\varphi \cdot 
\psi) $ for a linear functional $ \operatorname L$. 
 \end{remark}

\medskip

A central operator on  $ \ell^2(\mathbb N)$  the \emph{shift} operator 
$ \operatorname S (\alpha _0,\alpha _1,\ldots) \eqdef (0,\alpha _0,\alpha _1,\ldots)$. 
 Hankel operators $ \operatorname H$ are distinguished by their \emph{intertwining}  
 with the shift operator:
\begin{equation*}
\operatorname H \operatorname S = \operatorname S ^{*} \operatorname H\,.
\end{equation*} 
Proof is left for the reader.

The shift operator on Hardy space $ H ^2 (\mathbb D)$ is given by multiplication by $ z$, and Hankel operators  
on Hardy space enjoy the same intertwining with the shift operator.

While we have taken pains to outline these initial observations on the integers and the 
circle, there is an equivalent formulation on the real line.   To be specific, 
on $ L^2(\mathbb R )$, we have the Fourier transform 
\begin{equation*}
\widehat f(\xi )=\int f(x) \operatorname e ^{-i \xi x} \; dx  \,.
\end{equation*}
Define the orthogonal projections onto positive and negative frequencies 
\begin{equation*}
\operatorname P _{\pm} f(x) \eqdef \int _{\mathbb R _{\pm}} \widehat f(\xi ) 
\operatorname e ^{i \xi x} \; dx  \,.
\end{equation*}
Define Hardy spaces $ H ^2 _{\pm } (\mathbb R ) \eqdef \operatorname P _{\pm} L^2(\mathbb R ) $.
Functions $ f\in H^2_+(\mathbb R )$ admit an analytic extension to the upper half plane 
$ \mathbb C_+$. As in the case of the disk, it is convenient to refer to functions in 
$ H^2_+(\mathbb R )$ as \emph{analytic}.

A \emph{Hankel operator with symbol $ b$} is then a linear operator from $ H^2 _+(\mathbb R )$ 
to $ H^2 _{-}(\mathbb R )$ given by $ \operatorname H _{b} \varphi  \eqdef \operatorname P_- 
\operatorname M_b \overline{\varphi }$. 
This only depends on the analytic part of $ b$.

It will be convenient to consider some of our proofs in the setting of the real line. 
For, while it is equivalent to work on any of the three settings, the real line 
has a natural \emph{dilation structure} which simplifies certain aspects of the argument.

\begin{acknowledgment} These notes were prepared while in residence at the University 
of British Columbia, for a conference ``Harmonic Analysis at Sapporo, Japan''
 held in August 2005.
 \end{acknowledgment}

\section{Wavelets, $ \textup{BMO} (\mathbb R )$ and  Paraproducts } 

Ultimately, we are interested in characterizations of $ \textup{BMO}$.  This class of functions 
have delicate properties, sensitive to locations in both time and frequency.  \emph {Wavelets} 
turn out to be very useful in analyzing their behavior.    We recall some basic 
facts about two distinct classes of wavelets, namely the Haar and Meyer wavelets.

Throughout this paper,  $ \mathcal D$ denotes the dyadic grid.  Thus, 
\begin{equation}\label{e.dyadic}
\mathcal D \eqdef \{[j2^k,(j+1)2 ^{k})\mid j,k\in \mathbb Z \}\,.
\end{equation}

Define translation and dilation operators by 
\begin{align}\label{e.translate}
\operatorname {Tr} _{y}f(x) &\eqdef  f(x-y)\,,\qquad y\in \mathbb R \,,
\\ \label{e.dilate}
\operatorname {Dil} _{s} ^{p} f(x) &\eqdef s ^{-1/p} f(x/s)\,, 
\qquad 0<s,p<\infty \,,
\\ \label{e.transDil}
\operatorname {Dil} _{I} ^{p} f(x) &\eqdef \operatorname {Tr} _{c(I)} \operatorname {Dil}
_{\abs{ I}} ^{p}f(x)\,,\qquad \textup{ $ I$ is an interval}\,. 
\end{align}
In the second definition, $ s $ denotes the \emph {scale} of the dilation, and the normalization 
is chosen to preserve $ L^p (\mathbb R )$ norm.  In the last definition, we extend the definition of dilation 
to an interval, which incorporates a translation to the center of $ I$, denoted $ c(I)$, 
and a dilation by the scale of $ I$.

\subsection{Haar Wavelets} 

The notations $ h=h^0$ and $ h^1$ are reserved for the functions 
\begin{equation*}
h=h^0=-\mathbf 1 _{[0,1/2)}+\mathbf 1 _{[1/2,1]}\,,
\qquad 
h ^{1}=\mathbf 1 _{[0,1]}\,.
\end{equation*}
Here, the superscript $ {} ^0$ means that that the function has mean zero, while the 
superscript $ ^1$ means that the function does not have mean zero. 
We will use these two definitions in our discussion of paraproducts; the function 
$ h^0$ is used most of the time, and we will frequently suppress the superscript 
when using these functions. 

For any interval $ I$, we can define $ h ^{\epsilon }_I \eqdef \operatorname {Dil} _{I} ^{2} 
h ^{\epsilon }$. 
The Haar wavelets are then given by  $ \{h ^{0} _{I}\mid I\in \mathcal D\} $. 
These functions are an orthogonal basis on $ L^2 (\mathbb R )$, which 
extend to an unconditional basis on $ L^p (\mathbb R )$ for $ 1<p<\infty $.  
In particular, the Littlewood Paley inequalities become 

\begin{theorem}\label{t.lp-haar} We have the estimates 
\begin{equation*}
\norm f.p. \simeq \NOrm 
\Bigl[\sum _{I\in \mathcal D} \frac {\abs{ \ip f, h_I,} ^2} {\abs{ I}} \mathbf 1 _{I} 
\Bigr] ^{1/2} .p.\,, \qquad 1<p<\infty \,.
\end{equation*}
More generally, if $ \varphi $ is adapted to $ [0,1]$ and has mean zero, we have the 
square function below maps $ L^p$ into itself for all $ 1<p<\infty $.
\begin{equation}\label{e.LPS}
\operatorname S f \eqdef \Bigl[ \sum_{I\in \mathcal D} 
\frac {\abs{ \ip f, \varphi _I,} ^2} {\abs{ I}} \mathbf 1 _{I}\Bigr] ^{1/2} \,.
\end{equation}
\end{theorem}

One  natural extension of these estimates to the cases of $ p=1$ and $ p=\infty $ 
can be taken as the definitions of Hardy space $ H^1$ and $ \textup{BMO}$.

We should also mention that 
\begin{equation}\label{e.dyMax}
\operatorname M _{\textup{dy}} f(x)=\sup _{I\in \mathcal D} \frac 
 {\abs{ \ip f, h_I,} } {\sqrt{\abs{ I}}} \mathbf 1 _{I}
\end{equation}
is the dyadic maximal function.  It maps $ L^p$ into itself for all $ 1<p \le \infty $. 
The proof of this theorem can appeal to probabilistic methods.  Indeed, 
the expansion of a function in a Haar basis is a martingale.  
This fact lies behind the very successful application of Haar functions to 
establish a range of   deeper properties of singular integrals, including 
the Hilbert transform.  These properties include the   $ \textup{UMD}$ 
theory of Burkholder \cite{MR730072} and Bourgain \cite {MR727340} ;
matrix valued paraproducts, as discussed in articles by 
a range of authors  \cites {MR2000m:42016,MR1964822,MR98d:46039,MR2002m:47038,MR1880830}
and the Nazarov Treil Volberg extension of the Calder\'on Zygmund theory
\cites {MR1470373,MR1626935}, as well as their discussion of the Bellman function 
approach 
\cites {MR1428988,MR1945290}.

 More generally, the definition of the maximal function is 
\begin{equation}\label{e.maximal}
\operatorname M f(x) \eqdef  \sup _{t>0} (2t) ^{-1}\int _{-t} ^{t} f(x-y) dy
\end{equation}
where it is essential that we do not impose absolute values inside the integral. 
Define the 
Hilbert transform 
\begin{equation}\label{e.hilbert}
\operatorname H f(x) \eqdef -\operatorname P _{-}+\operatorname P _{+} f(x)
= 
\textup{p.v.}\tfrac1 \pi \int f(x-y) \frac {dy}y\,.
\end{equation}

\subsection{Meyer Wavelet}

Y.~Meyer \cite {MR98e:42001} found a Schwartz function $ w$, with 
\begin{equation}\label{e.meyer}
\textup{$ \widehat  w$  is supported on $  {2 \pi }\le \abs{ \xi  }\le 8 \pi $,}
\end{equation}
and the functions $ \{w_I\mid I\in \mathcal D\}$ form an orthonormal basis 
for $ L^2 (\mathbb R )$.  Here, we use the same notation as in the case of the 
Haar basis, $ w_I=\operatorname {Dil} _{I} ^{2}$. 

As with the Haar basis, these extend to an unconditional basis on $ L^p (\mathbb R )$, 
for $ 1<p<\infty $.  As concerns the Hilbert transform, observe that 
\begin{equation*}
\operatorname Hf = \sum_{I\in \mathcal D} \ip f,w_I,\, \operatorname H w_I,
\end{equation*}
and that the functions $ \operatorname H w_I= \operatorname {Dil} _{I} ^{2}(\operatorname H w)$
have the same decay and Fourier localization properties of the Meyer wavelet.  

Similarly, if $ f\in H^2 (\mathbb R )$, we have 
\begin{align*}
f= \operatorname P_+ f&=\sum _{I \in \mathcal D} \ip f, w_I, w_I
\\
&= \operatorname P_+ \Bigl[ \sum_{I\in \mathcal D} \ip  \operatorname P_+ f, w_I, w_I \Bigr]
\\
&= \sum_{I\in \mathcal D} \ip \operatorname P_+ f, w_I, \operatorname P_+w_I
\end{align*}
Therefore, $ \{ \operatorname P_+ w_I\mid I\in \mathcal D\}  $ is a basis for $ H^2 (\mathbb R )$. 

\subsection{Hardy Space $ H^1_+(\mathbb T )$ and $ \textup {BMO}$} 

For $ 1\le p\le \infty $, the Hardy space $ H^p (\mathbb T ) $ is the closure in 
$ L^p(\mathbb T )$ norm of the span of the exponentials $ \{z ^{n} \mid n\in \mathbb N\}$. 
A function $ f\in H ^p (\mathbb T )$ has an analytic extension $ F$ to the disk. 
Indeed, this extension is 
\begin{equation*}
F(z)=\sum _{n=0} ^{\infty } \widehat f(n) z ^{n}\,. 
\end{equation*}
and the $ H ^{p} (\mathbb T )$ norm can be taken to be 
\begin{equation*}
\norm f.H^p (\mathbb T ). \eqdef \sup _{0<r<1} \norm F(rz).L^p (\mathbb T ).
\end{equation*}

Concerning the space $ H^1 (\mathbb T )$, the following classical property 
is central to us. 

\begin{proposition}\label{p.factor}
Each function $ f\in H^1 (\mathbb T )$ is a product of functions $ f_1,f_2\in H^2 (\mathbb T )$, 
in particular, $ f_1$ and $ f_2$ can be chosen so that 
\begin{equation*}
\norm f.H^1.=\norm f_1.H^2.\norm f_2.H^2. 
\end{equation*}
\end{proposition}

\begin{remark}\label{r.factorfail} 
In the product setting to which we turn to next, this last property fails.  Indeed, 
part of the interest of our results is that while simple factorization will fail, 
 a notion of \emph {weak factorization} is in fact true. 
\end{remark}

\begin{remark}\label{r.rudin} The investigation of the failure of 
factorization is an intricate one.  In the product setting, Rudin \cite{MR0255841} proved the 
failure of factorization in the case of $ \mathbb D ^{d}$ for $ d\ge4$; 
Miles \cite{MR0374459} improved the result to $ d\ge3$; Rosay \cite{MR0377098} in the case of $ d=2$ showed that the 
set of functions in $ H^1(\mathbb D \otimes \mathbb D )$ which factor is of the 
first category.  This problem is also of interest in the Bergman space setting. 
See \cites{MR0427650,MR0338399,MR690048}. See \cite{MR690048} for information 
about this question in other spaces of analytic functions. 
\end{remark}

We are especially interested in the Hardy space $ H^1(\mathbb R )$.  
It is technically easier to discuss 
the \emph{real} Hardy space $ \operatorname {Re}(H^1)$ consisting of the real 
part of functions in $H^1$.

\begin{theorem}\label{t.H1} We have the equivalence of norms 
\begin{equation*}
\norm f. \operatorname {Re}(H^1).
\simeq \norm f.1.+\norm \operatorname H f.1.
\simeq \norm \operatorname Sf.1.
\simeq \norm \operatorname M  f.1.
\end{equation*}
Here $ \operatorname S$ is as in (\ref{e.LPS}).
\end{theorem}

Any standard reference in the subject will include a proof of this result. 
Historically, this kind of characterization was an essential precursor to the 
proof of Fefferman and Stein of $ H^1$ and $ \textup{BMO}$ duality.  

Indeed, one can use this characterization of $ H^1$ in terms of \emph {atoms} 
(which we don't define here) which lead immediately to 

\begin{theorem}[$ H^1 (\mathbb R )$---$ \textup{BMO} (\mathbb R )$ duality]
\label{t.dualit}  The dual of $ \operatorname {Re}( H^1 (\mathbb R) ) $ is $ \textup{BMO} (\mathbb R )$ with norm 
\begin{align*}
\norm f.\textup{BMO}(\mathbb R ).& \eqdef \sup _{ \textup{$ J$ is an interval}}
\Bigl[ \abs{ J} ^{-1} \sum_{\substack{I\in \mathcal D\\ I\subset J }} \abs{ \ip f,w_R,} ^2 \Bigr]
\\
&\simeq  \sup _{ \textup{$ J$ is an interval}} 
\Bigl[ \abs{ J} ^{-1} \int _{J} \abs{  f(x)-\mu _J} ^2 dx \Bigr]^{1/2} 
\end{align*}
where $ \mu _J= \abs{ J} ^{-1} \int _{J} f(y)\; dy$.
\end{theorem}

The second definition has the advantage of being intrinsic to $ f$, but has the disadvantage 
of not having a suitable generalization to higher parameters.  Thus, we have stressed 
the first definition in terms of wavelets.  There is nothing special about the 
Meyer wavelet appearing here. It can be replaced by any wavelet, including the Haar 
wavelet in this context.  

One nice feature of the Meyer wavelet, is that if we replace the functions 
$ w_R$ by their analytic projections, we obtain a completely analogous 
 definition of \emph {analytic} $ \textup{BMO}$, the dual to $ H ^{1} (\mathbb T )$.
.

\subsection{Paraproducts}

 \emph {Paraproduct} is the term used to refer to any of a wide variety of objects that 
 are a variant of a product of two functions.  A Hankel operator is 
 one such example, but our purpose in this section is to describe 
 a class of more naive examples.

 Consider the operation $ \operatorname M _{b} f$, where we take both $b  $  and $ f$ to have 
 finite Haar expansion 
 on the real line.  
 \begin{align*}
\operatorname M_b f = \sum_{I,J\in \mathcal D} \ip b,h_I, \ip f,h_J, h _I\cdot h_J
\end{align*}
Restrict the sum above to $ I\subsetneq J$, and observe that 
\begin{equation}  \label{e.para}
\begin{split}
\operatorname {Para} _{\textup{Haar}}(b,f)&=
\sum_{\substack{I,J\in \mathcal D\\ I\subsetneq J }} \ip b,h_I, \ip f,h_J, h _I\cdot h_J
\\
&= \sum_{I\in \mathcal D} \frac{\ip b,h_I, }{\sqrt {\abs{ I}} } \ip f,h^1_I, h_I 
\end{split}\end{equation}
Here, we are appealing to this property of Haar functions: 
\begin{equation*}
\ip h_J,h ^1 _I,= 
\begin{cases}
0 & I\cap J=\emptyset,\quad J\subset I,
\\
\frac{ h_J(c(I))}{\sqrt {\abs{ I} }} & I\subsetneq  J\,.
\end{cases}
\end{equation*}
The point of this is that while  the operator norm of $ \operatorname M_b$ is $ \norm b.\infty .$, 
while the norm of the operator in (\ref{e.para})
is in general somewhat smaller. 
$ \operatorname {Para} _{\textup{Haar}}(b,\cdot)$ is a paraproduct.

Paraproducts 
admit a more general definition, which is the point of this definition. 
For an interval $I $, we say that $\varphi $ is \emph {adapted to $I $} iff  $\norm \varphi.2.=1$ and 
 \begin{equation} \label{e.adapted} 
 \Abs{\operatorname D ^{n} \varphi(x)} {}\lesssim{}\abs{I}^{-n-\frac12} \Bigl( 1+\tfrac {\abs{x-c(I)}}{\abs{I}}\Bigr)^{-N}, 
 \qquad n=0,1. 
 \end{equation}
 Here, $c(I) $ denotes the center of $I $, and $N $ is a large  integer,
 whose exact value  need not concern us
  $\operatorname D$ denotes the derivative 
 operator. 
 We shall consistently work with functions which  have  $L^2$ norm at most one.  
 Some of these functions we will  also insist to have integral zero.

We take $ \{\varphi _ I^ \varepsilon \mid I\in \mathcal D\}$, $ \varepsilon \in \{0,1\}$
to be functions   adapted to $ I\in \mathcal D$.  The functions $ \varphi ^{0} _{I}$ are 
of mean zero. Then, set 
\begin{equation*}
\operatorname {Para}(b,f) \eqdef \sum_{I\in \mathcal D} 
 \frac{\ip b,\varphi ^0_I, }{\sqrt {\abs{ I}} } \ip f,\varphi ^1_I, \varphi ^0_I 
\end{equation*}

\begin{theorem}\label{t.para} 
For $ 1<p<\infty $, we have the inequality 
\begin{equation}\label{e.paranorm}
\norm \operatorname {Para} _{\textup{Haar}}(b,\cdot).p\to p. \lesssim \norm b. \textup{BMO}_{\textup{dy}}.,
\end{equation}
where the last norm is dyadic $ \textup{BMO}$ norm defined by 
\begin{equation}\label{e.dyBMO}
\norm b. \textup{BMO}_{\textup{dy}}. 
\eqdef 
\sup _{I \in \mathcal D} \Bigl[ \abs{ I} ^{-1} \sum_{J\subset I} 
\abs{ \ip b,h_I,} ^2\Bigr] ^{1/2} \,.
\end{equation}
For the operators $ \operatorname {Para}$, we have 
\begin{equation*}
\norm \operatorname {Para} (b,\cdot).p\to p. \lesssim \norm b. \textup{BMO}.,
\end{equation*}
\end{theorem}

It is essential in this formulation that we do not let the function 
$ h ^{1} _{I} $, which does not have mean zero, fall on the $ \textup{BMO}$ function $ b$.
This point of view is very helpful in obtaining upper bounds for other commutators, 
Hankel operators and related objects.

\bigskip 

It is useful to observe that paraproducts can arise in a wide variety of forms, 
in particular, the classical approach of Coffman and Meyer relies upon the 
`$\operatorname P_t $--$ \operatorname Q_t$' formalism.  A form useful to us 
is as follows.  For the Meyer wavelet $ w$, with it's antianalytic and analytic 
parts, respectively $ u$ and $ v$, let us set 
\begin{align}\label{e.deltaU}
\Delta \operatorname U_j &\eqdef \sum_{\substack{I\in \mathcal D\\ \abs{ I}= 2 ^{j} }}
u _{I}\otimes u_I 
\\ \label{e.Uj}
\operatorname U_j & \eqdef \sum _{k\ge j} \Delta \operatorname U_k 
\end{align}
The following Theorem concerns a paraproduct which is quite close to being 
a Hankel operator;  it's bounds are a  corollary to the previous Theorem. 

\begin{theorem}\label{t.U} We have the estimate 
\begin{equation}\label{e.U}
\NOrm \sum_{j\in \mathbb Z } (\Delta \operatorname U_j b)\cdot \overline{(\operatorname U_j) \varphi } 
.2. \lesssim \norm b.\textup{BMO}(\mathbb R ).\norm \varphi .2.
\end{equation}
\end{theorem}

\begin{proof}

Let us first consider terms like 
\begin{equation*}
\sum _{j \in \mathbb Z } 
\Delta \operatorname U_j  \cdot  \overline{\Delta \operatorname U _{j+k}}\,,
\qquad 
0\le k\le8\,.
\end{equation*}
We take e.g.~$ \abs{ k}\le8$ due to the compact frequency support of the 
Meyer wavelets, as will be come clear momentarily.  The assertion 
is that each of these is a bounded operator, provided $ b\in \textup{BMO}(\mathbb R )$. 

We control this expression in a brute force method, which we will appeal to 
twice. 
Fix $k$ and 
 a map $ \pi \mid \mathcal D
\longrightarrow \mathcal D$ so that $ \abs{ \pi (I)}=2^k I$ 
and $(A-1) \abs{ I}\le{} \operatorname {dist}(I,J)<A\abs{ I}$, for some fixed 
integer $ A$.  

Write   
\begin{gather*}
\psi _{I} \eqdef \sqrt {\abs{ I}} \, u_I \cdot \overline{u _{\pi (I)}}
\end{gather*}
Observe that $ A ^{100} \psi _{I}$ is adapted to $ I$ with constant 
independent of $ A$ (and we certainly do not 
assert that it has integral zero!). This is the only observation 
we need to make to conclude that 
\begin{equation*}
\NOrm \sum _{I\in \mathcal D} \frac {\ip b, u_I,} {\sqrt {\abs{ I}}} 
\overline{\ip \varphi , u_ {\pi (I)} ,}  \psi _{I}
.2. \lesssim A ^{-100} \norm b. \textup{BMO}(\mathbb R ). 
\norm \varphi .2.\,.
\end{equation*}
This is summed over the different values of $ A$ and $ \pi $ to conclude the 
estimate 
\begin{equation*}
\NOrm 
\sum _{j \in \mathbb Z } 
\Delta \operatorname U_j b  \cdot  \overline{\Delta \operatorname U _{j+k}\varphi }\,,
.2. \lesssim 
\norm b. \textup{BMO}(\mathbb R ). 
\norm \varphi .2.\,, \qquad 0<k<8\,.
\end{equation*}

Associated with the Meyer wavelet is a `father wavelet,' a function  $ W$ of 
non zero mean, for which $ w= W- \operatorname {Dil}_ {1/2} ^{1}W$. Using this,   
we see that 
\begin{equation*}
\operatorname U_j f=\sum _{\abs{ I}=2^j} \ip f, W_I,\, W_I\,. 
\end{equation*}
We then have 
\begin{equation*}
\sum _{j\in \mathbb Z } (\Delta \operatorname U_j b)
\cdot \overline{ \operatorname U _{j+9} \varphi }
= \sum_{\substack{I,J\in \mathcal D\\ 2^9 \abs{ I}\le \abs{ J} } }
\frac {\ip b, u_I,} {\sqrt {\abs{ I}}} \,
\overline{ \ip \operatorname P_+ \varphi , W_J, }\,\sqrt {I} \, u_I \cdot W_J\,.
\end{equation*}
Observe that in contrast to the previous case, 
we are taking inner products $ \ip  \operatorname P_+ f,W_J,$ and $ W_J$ will have 
a non zero mean. 

Nevertheless, this last sum can be controlled by the brute force method used above, with this observation. 
Take two dyadic  intervals $ I$ and $ J$ with $ 2^9 \abs{ I}=\abs{ J}$ and 
$ (A-1) \abs{ I}\le{} \operatorname {dist} (I,J)\le{} A \abs{ I}$, for integer $ A$. 
Then, 
\begin{equation*}
A ^{100}\sqrt {I} \, u_I \cdot W_J 
\end{equation*}
is adapted to $ I$ with constant independent of $ A$, \emph {and has integral zero.}
The reason for this  stems from (\ref{e.meyer}). The Fourier transform of the product is  supported in the 
convolution of the supports of the Fourier transforms of the two functions.  
Since $ I$ has the smaller scale, the Fourier transform of $ u_I$ is supported a 
very great distance from the origin, hence the Fourier transform is zero at the origin. 
This completes the proof. 
\end{proof}

\section{The Nehari Theorem on the  Disk} 

The classical result that we are interested in is: 

\begin{nehari}[\cite {nehari}] \label{t.nehari}
The Hankel operator $ \operatorname H_b$ from  $ H^2_+(\mathbb T )$ to $ H^2_ +(\mathbb T )$ 
iff there is a bounded function $ \beta $ with $ P_+b=P_+\beta $.  Moreover, 
\begin{equation*}
\norm \operatorname H_b..=\inf _{\beta \mid \operatorname P_+\beta =\operatorname P_+ b} 
\norm \beta .\infty .
\end{equation*}
\end{nehari}

There are three proofs of this fact in the literature. 
In the new results, we will need to rely upon methods from two of these methods.

\subsubsection*{Factorization}

Given a bounded Hankel operator $ \operatorname H _{b}$, we want to show that we 
can construct a bounded function $ \beta $ so that the analytic part of $ b$ and 
$ \beta $ agree. 

This proof is the one found by Nehari \cite {nehari}.  We begin with a basic computation of the 
norm of the Hankel operator $ \operatorname H _{b}$: 
\begin{equation}\label{e.dualNorm}
\begin{split}
\norm \operatorname H _{b}..&=\sup _{\norm \varphi  . H^2_+ (\mathbb T ).=1}
\sup _{\norm \psi . H^2_+ (\mathbb T ).=1}
\int  \operatorname H_b \psi \cdot \overline{ \varphi} \; dx  
\\
&=\sup _{\norm \varphi  . H^2_+ (\mathbb T ).=1}
\sup _{\norm \psi . H^2_+ (\mathbb T ).=1} 
\int \operatorname P_+ \operatorname M_b \overline{\psi}  \cdot \overline \varphi \; dx 
\\
&=\sup _{\norm \varphi  . H^2_+ (\mathbb T ).=1}
\sup _{\norm \psi . H^2_+ (\mathbb T ).=1}\int 
(\operatorname P_+ b)\overline{\psi \cdot \varphi}  \; dx 
\\
&= \sup _{\norm \varphi  . H^2_+ (\mathbb T ).=1}
\sup _{\norm \psi . H^2_+ (\mathbb T ).=1} \ip (\operatorname P_+ b), \psi \cdot \varphi ,
\end{split}
\end{equation}
But, the $ H^1 (\mathbb T )=H^2 (\mathbb T )\cdot H^2 (\mathbb T )$, 
as we recalled in Proposition~\ref{p.factor}. 
We read from the equality above that 
the analytic part of $ b$ defines a bounded linear functional 
on $ H^1 (\mathbb T )$ a subspace of $ L^1 (\mathbb T )$.  

The Hahn Banach Theorem applies, giving us an extension of this linear functional 
to all of $ L^1$, with  the same norm.  But a linear function on $ L^1$ is a bounded 
function, hence we have constructed a bounded function $ \beta $ with the same analytic 
part as $ b$.

\subsubsection*{Duality}

In this proof, the $ H^1$---$ \textup{BMO}$ duality is decisive. 
The calculation (\ref{e.dualNorm}) shows that $\operatorname  P _{+} b$ is a bounded 
linear functional on $ H^1$.  Therefore, we have 
\begin{equation*}
\norm \operatorname H_b .. \simeq \norm \operatorname P_+ b . \textup{BMO}. 
\end{equation*}
(This is not equality, since we are not choosing to adopt a canonical norm for $ \textup{BMO}$.) 
In addition, we have $ \textup{BMO}=L^\infty + \operatorname H L^\infty $, where 
$ \operatorname H$ is the Hilbert transform.  Therefore, we can select $ \beta \in L^\infty $
which has the same analytic part as $ b$.

\begin{remark}\label{r.duality} 
Historically, this proof came last;  it depends critically upon the 
Fefferman Stein $ H^1$--$ \textup{BMO}$ duality, which was not established until the 
1970's, see \cite {MR56:6263}.
\end{remark}

\subsubsection*{The AAK Method}

Adamjan, Arov and Krein \cite {MR0298453} invented a method based upon 
an dilation\footnote {This dilation property is distinct from the dilation property 
of the real line used in other parts of this paper.}
property of operators on Hilbert space.  This method avoids the 
finer aspects of Hardy spaces.

We are given a Hankel matrix which is bounded from  $ \ell^2 (\mathbb N)$
to itself, and we seek to extend it to a bounded matrix on $ \ell^2 (\mathbb Z )$, 
with the same norm.  This is an inductive process, with the first step being 
that we seek to add, say, a row to the `top' of the matrix: 
\begin{equation*}
\begin{bmatrix} 
a_0 & a_1 &  a_2 & \dots 
\\
 a_1 &  a_2 & a_3 & \dots 
\\
 a_2 & a_3 &  a_4 & \dots 
\\
\vdots & \vdots & \vdots &  \ddots
\end{bmatrix}
\longrightarrow 
\begin{bmatrix} 
*   & a_0 & a_1 &    \dots  
 \\
 a_0 & a_1 &  a_2 & \dots 
\\
 a_1 &  a_2 & a_3 & \dots 
\\
 a_2 & a_3 &  a_4 & \dots 
\\
\vdots & \vdots & \vdots &  \ddots
\end{bmatrix}
\end{equation*}
Namely, we seek to choose a value of $ *$ to put on the upper left hand coordinate so that 
the two matrices have the same norm.

This in fact can be done, and 
leads to the following Proposition.

\begin{proposition}\label{p.aak}
Consider two Hilbert spaces $ \mathcal G$ and $ \mathcal H$, and consider 
linear operators from $ \mathcal G \oplus \mathcal H$ into itself of the form 
\begin{equation*}
\operatorname U=
\begin{bmatrix} 
\operatorname X & \operatorname C \\ \operatorname  A &  \operatorname B
\end{bmatrix}
\end{equation*}
where $ \operatorname X\mid \mathcal G\to \mathcal G $; $ \operatorname A \mid 
\mathcal G\to \mathcal H$; $ \operatorname B\mid \mathcal H\to \mathcal G $; 
and $ \operatorname C\mid \mathcal H\to \mathcal G $.  We presume that $ \operatorname A$, 
$ \operatorname B$, $ \operatorname C$ are prescribed in advance. 
Then we can select $ \operatorname  X$ so that 
\begin{equation*}
\norm \operatorname U..=\max\{ \norm A..\,,\ \norm B..\,,\ \norm C..\}. 
\end{equation*}
\end{proposition}

Given a Hankel matrix $ \operatorname H=\{a_{j+k}\mid j,k\in \mathbb N\}$, 
we apply the proposition above with 
\begin{equation*}
\operatorname A=\begin{bmatrix}
a_0 \\ a_1 \\ a_2 \\ \vdots \end{bmatrix},\qquad 
\operatorname B=
\begin{bmatrix}  a_1 & a_2 & a_3 & \dots 
\\
a_2& a_3 & \dots 
\\
\vdots & \vdots & \ddots 
\end{bmatrix}
,
\qquad 
\operatorname C=\begin{bmatrix} a_0 & a_1 & a_2 & \dots 
\end{bmatrix}
\end{equation*}
By the proposition, we can choose $ a _{-1}$ so that the norm of $ \operatorname H$ is 
the norm of 
\begin{equation*}
\begin{bmatrix}a _{-1} & a _{0} & a _{1} & \dots 
\\ 
  a _{0} & a _{1} & a_2 & \dots
 \\
 a _{1} & a_2 &  a_3 & \dots
 \\
 \vdots & \vdots & \vdots & \ddots 
\end{bmatrix}
\end{equation*}

By induction, we can extend the Hankel operator to $ \ell^2(\mathbb Z )$ to itself, as 
a bounded operator with the same norm.  The conclusion of Nehari's theorem then follows.

\begin{remark}\label{r.aak} This method has found many deep extensions to Hankel 
matrices whose entries are themselves operators.  We refer the reader to 
Nikolski \cite {MR1864396}, as well as Nikolski \cite{MR1892647} and Peller \cite{MR1949210} 
for very interesting discussions of the method of Adamyan, Arov and Krein. 
This is relevant for us, as  in the extension to the 
product setting, we consider Hankel matrices whose entries are also Hankelian.  
Cotlar and Sadosky have studied extensions of this method to the polydisk in 
a sequence of papers \cites{MR1066468
,
MR1067439
,
MR1203463
,
MR1233667
,
MR1300214
,
MR1415032}
\end{remark}

\subsection{Commutators}

Commutators are very useful in measuring, in a quantitative way, the distance 
from being abelian.  They are relevant for us as the 
Nehari theorem has an equivalent formulation in terms of 
commutators of multiplication operators and the Hilbert transform.  

\begin{theorem}\label{t.comm} We have the equivalence 
\begin{equation*}
\norm [\operatorname M _{b}, \operatorname H].2\to 2.\simeq \norm b. \textup{BMO}. 
\end{equation*}
Here, $  \textup{BMO}=(\textup{Re}H^1(\mathbb T ))^\ast$ is real $\textup{BMO} $.
\end{theorem}

The proof in this circumstance is immediate: Observe that  $ 
\operatorname P _{\mp}[\operatorname M _{b},\operatorname H] P _{\pm } $ is itself a Hankel 
operator, or a conjugate of a Hankel operator.  Specifically, 
\begin{align*}
\operatorname P_+ [\operatorname M _{b},\operatorname H] \operatorname   P_+&=0
\,, \qquad &
\operatorname P_- [\operatorname M _{b},\operatorname H] \operatorname   P_-&=0 \,,
\\
\operatorname P_+ [\operatorname M _{b},\operatorname H] \operatorname   P_-&=-\operatorname P_+
\operatorname M _{b} \operatorname P_-\,, \qquad &
\operatorname P_- [\operatorname M _{b},\operatorname H] \operatorname   P_+&=\operatorname P_-
\operatorname M _{b} \operatorname P_+\,.
\end{align*}
and the last two operators are orthogonal, and Hankel operators. 

The Theorem admits many extensions.  For instance, we continue to have the equivalence 
\begin{equation*}
\norm [\operatorname M _{b}, \operatorname H].p\to p.\simeq \norm b. \textup{BMO}. 
\,, \qquad 1<p<\infty \,.
\end{equation*}
Indeed, assuming the commutator with symbol $ b$ is bounded on $ L^p$, the 
same is true on the dual index $ L ^{p'}$.  But then, by interpolation, the commutator 
is bounded on $ L^2$, and we can appeal to the Theorem above to deduce that the symbol 
is in $ \textup{BMO}$.  

The upper bound on the $ L^p$ norm of the commutator follows by considering the 
Hankel operator $ \operatorname P_-
\operatorname M _{b} \operatorname P_+$. 
Using the calculation in (\ref{e.dualNorm}), we can see that 
\begin{equation*}
\norm \operatorname P_-
\operatorname M _{b} \operatorname P_+.p\to p.=
\sup _{\varphi \in H^p} \sup _{\psi \in H^{p'}} \int (\operatorname P _{-} b)\cdot 
\varphi \psi \; dx  
\end{equation*}
But, $ H^1=H^p\cdot H ^{p'}$, so this last quantity is $ \norm \operatorname P_- 
b. \textup{BMO}.$.

\medskip 

One direction in which this result extends is for the commutator to characterize 
a broad array of function spaces. The genesis of this theme is the 
very interesting article of Coifman, Rochberg and Weiss \cite{MR54:843}, 
which consider the instance of commutators of $ \operatorname M _{b}$ and 
Reisz transforms.  

Subsequently, it turns out that one has the equivalence 
\begin{equation*}
\norm [ \operatorname M _{b}, \operatorname H].X\to Y. \simeq \norm b.Z.
\end{equation*}
for a range of spaces $ X$, $ Y$, and function spaces $ Z$.  Various Lipschitz 
classes can be characterized this way; further generalizations can be stated in terms 
of various Besov and Treibel Lizorkin spaces. 
Moreover, the Hilbert transform can be 
replaced by other operators, such as the fractional integral operators. 
There is a significant literature here, of which we cite 
Chanillo \cite {MR84j:42027}, Cruz-Uribe and Fiorenza \cite {MR2004a:42010}, 
as well as the references in the article of the author \cite {math.CA/0502336} 
which is a first step in extending some of these results to 
higher parameter settings.

\medskip 

Another direction is to abandon the chance of characterizing function spaces, replacing 
the Hilbert transform by, say, a Calder\'on Zygmund operator $ \operatorname T$.
The method of choice 
in such generalizations is the use of the sharp function:
\begin{equation*}
([ \operatorname M _{b},\operatorname T]f) ^{\sharp } \lesssim \operatorname M f,
\end{equation*}
where on the right we intend $ \operatorname M$ to be an appropriate maximal function. 
This method was (to the best of my knowledge) first used  on this problem by 
Coifman, Rochberg and Weiss \cite {MR54:843}, and since then has been used by a wide variety of 
authors.

This method has difficulties in being generalized to the product setting; we ask the 
reader's patience for not defining precisely what the sharp function is.\footnote{
A paper of R.~Fefferman \cite {MR90e:42030} indicates a 
certain extension sharp function  to the two parameter product setting. 
The difficulty of a similar extension to a three parameter setting centers around the 
tenuous relationship between rectangular $ \textup{BMO}$ and $ \textup{BMO}$ in three 
and higher parameters. See however \cite{MR1767858}.}

There is an alternate method, which is of interest as it highlights the role of 
paraproducts in these questions, and permits a generalization to the higher parameter 
setting. See the author's paper \cite {math.CA/0502336}.

\subsection{Paraproducts and Commutators}

We describe how to write the commutator $ [\operatorname M _{b},\operatorname H]$ 
as a sum of two paraproducts.  From this, appropriate upper bounds on the
norms of the commutator can be given.  In order to keep the exposition as simple 
as possible, we will rely on particular properties of the Hilbert transform.  Yet, the 
method is flexible, and can apply to a wide variety of operators.  

We return to the Haar basis on the real line.  For convenience, set 
\begin{equation*}
g_I=-h _{I _{\textup{left}}}+h _{I _{\textup{right}}},\qquad I\in \mathcal D\,.
\end{equation*}
Here, we are using obvious notation to refer to the left and right halves of 
the dyadic interval $ I$, which again are dyadic.  
And define an operator by $ \operatorname G f=\sum _{I\in \mathcal D} \ip f,h_I,g_I$.

Now, it is well known that the Hilbert transform is nearly diagonalized in 
a wavelet basis. Yet, since the Hilbert transform has odd kernel, it is is not appropriate 
to decompose with a `even' choice of basis. Clearly, $ g_I $ is an `odd' version of the 
Haar function $ h_I$, so $ \operatorname G$ is a bit like the Hilbert transform. 
Of course it lacks translation and dilation invariance. 
It is a very nice observation of S.~Petermichl \cite {MR2000m:42016} (also see \cite {MR1964822}) 
that  these are the only properties missing.  Namely, we have 

\begin{proposition}[Petermichl \cite{MR2000m:42016}]\label{p.petermichl}
The operator below is a 
non zero multiple of the Hilbert transform. 
\begin{equation} \label{e.A}
 \lim _{Y\to \infty }  \int_{0} ^{Y} \int _{1} ^{2}  
\operatorname {Tr} _{-y} \operatorname {Dil} _{1/s} ^{(2)} 
\operatorname G 
\operatorname {Dil} _{s} ^{(2)}\operatorname {Tr} _{y} \; \frac{dy}y\frac {ds}s  
\end{equation}
\end{proposition}

\begin{proof}
Observe that the limit 
\begin{equation*}
\lim _{Y\to \infty }  \int_{0} ^{Y} \int _{1} ^{2}  
\operatorname {Tr} _{-y} \operatorname {Dil} _{1/s} ^{(2)} 
\operatorname G  \varphi 
\operatorname {Dil} _{s} ^{(2)}\operatorname {Tr} _{y} \; \frac{dy}y\frac {ds}s  
=\Lambda \varphi 
\end{equation*}
exists for Schwartz functions $ \varphi $.  Moreover, as $ G$ is a bounded operator 
on $ L^2$, we conclude that $ \Lambda $ is also bounded on $ L^2$.  

Due to the limiting procedure, one sees that $ \Lambda $ is a translation invariant 
operator.  The average over dilations is taken with respect to Haar measure 
for the dilation group, hence $ \Lambda $ is also invariant with respect to 
dilations.   It is a classical fact that a bounded linear operator, invariant 
with respect to translations and dilations is a linear combination of the 
Identity operator and the Hilbert transform.   

Clearly, $ \Lambda \mathbf 1 _{}=0$.  That is, $ \Lambda $ is 
a multiple of the Hilbert transform. And so we should check that it is a not 
the zero operator.   But one can check directly that  $ \Lambda $ applied 
to the Dirac measure at the origin is 
\begin{equation*}
\Lambda \delta _0 (x) \simeq x ^{-1}.
\end{equation*} 
Hence it is a multiple of the Hilbert transform. 
\end{proof}

\begin{proposition}\label{p.comm->para} 
The commutator $  [ \operatorname M _{b}, \operatorname G]$ can be written as a 
linear combination of paraproducts, or paraproducts composed with $ \operatorname G$.  
In particular, the commutator is bounded on $ L^2$, when $ b\in \textup{BMO}$. 
\end{proposition}

\begin{proof}
We use the notation $ \psi \otimes \varphi $ to denote the rank one linear operator 
\begin{equation*}
\psi \otimes \varphi (f)= \psi \ip \varphi ,f,\,. 
\end{equation*}

We will expand the symbol $ b$ in the Haar basis.  $ \operatorname G$ is an explicit 
sum over rank one operators as above, and we will make an computation of commutators 
for Haar functions.  As such, it is convenient to split the operator $ \operatorname G$ 
into $ \operatorname G _{\textup{left}}$ and $ \operatorname G _{\textup{right}}$, 
where we define 
\begin{equation*}
 \operatorname G _{\textup{left}} \eqdef 
 \sum _{I} h _{I _{\textup{left}}} \otimes h_J 
\end{equation*}
with a similar definition for $ \operatorname G _{\textup{right}}$.  Below 
we will only consider the `left' version.

We are lead to expand the commutators 
\begin{align} \nonumber 
[ \operatorname M _{h_I}, h _{J _{\textup{left}}}\otimes h_J  ]
& = 
(h_I h _{J _{\textup{left}}})\otimes h _{J} 
	- h _{J _{\textup{left}}}\otimes (h_I h_J)
\\  \label{e.commCases}
&= \abs{ J} ^{-1/2} \begin{cases}
0 & I\cap J=\emptyset,\quad J\subsetneq I 
\\
\sqrt 2\bigl[{h _{J _{\textup{left}}} ^{1}}  
\otimes h_J  \pm h _{J _{\textup{left}}}
\otimes  {h _{J _{\textup{left}}} }  \bigr]
& I= J _{\textup{left}}
\\
  {h _{J _{\textup{left}}} }  
\otimes h _{J _{\textup{right}}} 
& I= J _{\textup{right}}
\\
-\sqrt 2\; {h _{J _{\textup{left}}} }  
\otimes h_J  
-h _{J _{\textup{left}}} \otimes 
{h _{J } ^{1}}  
& I=J
\\
\sqrt 2 \; h_I \otimes  {h _{J }}  
- {h _{J _{\textup{left}}} }   \otimes h_I
& I\subsetneq J _{\textup{left}}\,.
\end{cases}
\end{align}
In this computation, we note that there are two conditions that lead 
to the commutator being zero.  The first is a trivial localization condition, 
$ I\cap J=\emptyset$.  The second condition, $ J\subsetneq I$,
is an essential cancellation  condition coming from commutator.  

All other terms lead to a paraproduct term, although some of these paraproducts 
are trivial, in that all relevant functions have a zero. Apply this computation to the
 commutator in question, expanding as 
 \begin{align*}
[ \operatorname M _{b}, \operatorname G _{\textup{left}}]=
\sum _{I,J\in \mathcal D} \ip b, h_I, 
[ \operatorname M _{h_I}, h _{J _{\textup{left}}}\otimes h_J  ]
\end{align*}

For instance, from the case of $ I= J$, we get  
\begin{equation*}
c\sum _{J} \frac{\ip b,h_J,}{\sqrt {\abs{ J }}} \;
h _{J _{\textup{left}}} \otimes 
{h _{J _{\textup{left}}} ^{1}} 
\end{equation*}
which is a paraproduct operator with symbol $ b$. 
Notice that the `$ 1$' falls on the right side of the tensor product of Haar functions. 
In considering the case $ I={J _{\textup{left}}}$, we get a paraproduct that is dual 
to the one above, namely 
\begin{equation*}
\sum _{J} \frac{\ip b, h _{J _{\textup{left}}},} {\sqrt {\abs{ J }}}  \;
{h _{J _{\textup{left}}} ^{1}}
\otimes h_J \,.
\end{equation*}
The other term that arises 
from the case $ I={J _{\textup{left}}}$ is less singular:
\begin{equation*}
\sum _{J} \frac{\ip b, h _{J _{\textup{left}}}, } {\sqrt {\abs{ J }}} \;
h _{J _{\textup{left}}} 
\otimes h_J 
\end{equation*}
Here, all functions are Haar functions, that is they have zeros.  This case is 
easier to control. 

Let us consider the case of $ I\subsetneq J _{\textup{left}}$.  
Observe that we have 
\begin{equation*}
 \operatorname G _{\textup{left}} ^{\ast} h ^{1} _{I}=-
 \sum _{I\subsetneq J _{\textup{left}}}  \sqrt{ \tfrac 
 {\abs{ I}} {\abs{ J}} } \; h_J
\end{equation*}
Keeping this in mind, we see that 
\begin{align*}
 \sum _{I\subsetneq J _{\textup{left}}} 
\frac{\ip b, h_I,}{\sqrt {\abs{ J }}}   \;
h_I \otimes  {h _{J }}      
& = -
\sum _{I } 
\frac{\ip b, h_I,}{\sqrt {\abs{ I }}}   \;
h_I    \otimes ( G^\ast h_I^1)  
\end{align*}
That is, we have a composition of $ \operatorname G ^{\ast }$ and a paraproduct. 
For the other term associated with this case we have 
\begin{align*}
\sum _{I\subsetneq J _{\textup{left}}} 
\ip b, h_I,  
 \frac{h _{J _{\textup{left}}} } {\sqrt {\abs{ J }}} \otimes h_I
\end{align*}
This is dual to the previous case.
 Our proof is complete. 

\end{proof}

\begin{remark}\label{r.fracInt} 
A simpler exposition of this approach can be had for commutators of multiplication 
operators and fractional integral operators.  See Lacey \cite{math.CA/0502336}.  This provides an alternate proof 
of a result of Chanillo.  This result proves another characterization of $ \textup{BMO}$ 
in terms of a commutator.  
\end{remark}

%
%
%
%
%
%


\section{Aspects of Product Hardy Theory} 

We describe the elements of product Hardy space theory, 
as developed by S.-Y.~Chang and R.~Fefferman \cites{MR86g:42038,MR82a:32009,
MR90e:42030,MR86f:32004,MR81c:32016} as well as Journ\'e \cites{MR87g:42028
,
MR88d:42028}.  
By this, we mean 
the Hardy spaces associated with domains like $ \mathbb D \otimes \mathbb D $, 
with boundary $ \mathbb T \otimes \mathbb T $. In particular, the boundary is flat, 
and while we work with several variables, we are very far from the pseudoconvex case.

 We view $ \mathbb R ^{d}$ as a tensor 
product of one dimensional spaces.  In particular, previously, we used the 
splitting of $ L^2 (\mathbb R )=H^2 (\mathbb 
C_+ )\oplus H^2 _{-} (\mathbb C_+ )$.  
This leads to a decomposition of 
\begin{equation*}
 L^2 (\mathbb R ^{d}) =\bigotimes _{j=1} ^{d} L ^2 (\mathbb R ) =
 \bigotimes _{j=1} ^{d} H^2 (\mathbb C_+ )\oplus H^2 _{-} (\mathbb C_+ )\,,
\end{equation*}
into $ 2 ^{d}$ components.  

To describe them, let us set $ \operatorname P _{\pm,j}$ to be the 
one dimensional Fourier projection operator $ \operatorname P _{\pm}$
acting on the $ j$th coordinate.  For $ \sigma \in \{-,+\} ^{d}$, set 
\begin{equation*}
\operatorname P _{\sigma }=\bigotimes _{j=1} ^{d} \operatorname P _{\sigma (j),j}
\end{equation*}
Likewise, we set $ H^2 _{\sigma } (\mathbb C_+^{d})$ to be the range of the 
orthogonal projection $ \operatorname P _{\sigma }$.  We then 
have 
\begin{equation*}
L^2 (\mathbb R ^{d})=\bigoplus _{\sigma \in \{+,-\} ^{d}}
H _{\sigma } ^{2} (\mathbb C_+^{d}) 
\end{equation*}

Among these $ 2^d$ Hardy spaces, we distinguish $ H _{\oplus } ^2 (\mathbb C_+^{d})$ 
in which $ \sigma \equiv +$, and likewise for $ H _{\ominus} ^2 (\mathbb C_+^{d})$. 
The corresponding orthogonal projections are $ \operatorname P _{\oplus}$ and $ 
\operatorname P _{\ominus}$.

Functions $ f$ in this space are defined on $ \mathbb R ^{d}$.
 $\mathbb R^d$ is viewed as the boundary of the `upper half space'
 \begin{equation*}
\mathbb C_+^d=\prod_{j=1}^d \{z\in\mathbb C\mid \operatorname{Re}(z)>0\}
\end{equation*}
And we require that there is a function $F\,:\,\mathbb C_+^d\longrightarrow\mathbb C$ 
that is holomorphic in each variable separately, and 
 \begin{equation*}
f(x)=\lim_{\norm y..\to0}f(x_1+iy_1,\ldots,x_d+iy_d).
\end{equation*}	
The norm of $f$ is taken to be 
 \begin{equation*}
\norm f. \operatorname H^p _{\oplus} (\mathbb C_+ ^{d}).=\lim_{y_1\downarrow 0}\cdots\lim_{y_d\downarrow 0}\norm F(x_1+y_1,\ldots,x_d+y_d).L^1(\mathbb R^d).
\end{equation*}

\begin{remark}\label{r.H(Rd)}
The (real) Hardy space $ H^1 (\mathbb R ^{d})$ typically denotes the class of functions 
with the norm 
\begin{equation*}
\norm f.1.+\sum _{j=1} ^{d} \norm \operatorname R_j f.1.
\end{equation*}
where $ \operatorname R_j$ denote the Reisz transforms. This space is 
invariant under the one parameter family of isotropic dilations, while   $ H^1(\mathbb C_+^d) $ 
is invariant under dilations of each coordinate separately. That is, it is invariant 
under a $ d$ parameter family of dilations.  That is why we refer to  `multiparameter' 
theory, or `$ d$ parameters.'
\end{remark}

As before, the real $ H^1$,  $\operatorname {Re} H^1(\mathbb C_+^d) $ has a variety of equivalent norms, 
in terms of square functions, maximal functions and Hilbert transforms.  
For our discussion of paraproducts, it is appropriate to make some definitions of translation 
and dilation operators which extend the definitions in (\ref{e.translate})---(\ref{e.transDil}). 
(Indeed, here we are adopting broader notation than we really need, in anticipation 
of a discussion of multiparameter paraproducts.) 
Define 
\begin{align}\label{e.translated}
\operatorname {Tr} _{y}f(x-y) &\eqdef f(x-y),\qquad y\in \mathbb R ^{d}\,,
\\ \label{e.dilated}
\operatorname {Dil} _{t_1,\dotsc,t_d}^p f(x_1,\dotsc,x_d) 
&\eqdef (t_1\cdots t_d) ^{-1/p} f(x_1/t_1,\dotsc,x_d/t_d)\,,\qquad t_1,\dotsc,t_d>0
\\ \label{e.transDild}
\operatorname {Dil} _{R} ^{p} &\eqdef \operatorname {Tr} _{c(R)} \operatorname {Dil} 
_{\abs{ R_1},\dotsc , \abs{ R_d}} ^{p}\,.
\end{align}
In the last definition $ R=R_1\times\cdots\times R_d$ is a rectangle, and the dilation 
incorporates the locations and scales associated with $ R$.  $ c (R)$ is the center of $ R$.

Let $ \mathcal D ^{d}=\mathcal D\times \cdots\times \mathcal D$ denote the $ d$ fold 
product of the dyadic intervals.  These are the dyadic rectangles in $ \mathbb R ^{d}$. 
For a non negative bump function $ \varphi ^{1}$ with $ \int \varphi ^{1}\; dx=1$, 
define the (strong) maximal function by 
\begin{equation*}
\operatorname M\cdots \operatorname M f(x)=\sup _{R\in \mathcal D ^{d}} 
\operatorname {Dil} _{R} ^{2} \varphi^1 (x) \ip f,\operatorname {Dil} _{R} ^{2}\varphi ^{1},  
\end{equation*}
We use the superscript on $\varphi ^{1}  $ to indicate that it has a non zero integral.

Fix a bump function $ \varphi ^{0}$ so that 
\begin{equation*}
\varphi ^{0}(x_1,\dotsc,x_d)=\prod _{j=1} ^{d}\varphi (x_j) 
\end{equation*}
where $ \int _{\mathbb R } \varphi \; dx=0$.  Then set an analog of the Littlewood 
Paley square function to be 
\begin{equation*}
\operatorname S\cdots \operatorname S f(x)=
\Bigl[\sum _{R\in \mathcal D ^{d}} 
[\operatorname {Dil} _{R} ^{2} \varphi^0 (x)] ^{2}
\abs{ \ip f,\operatorname {Dil} _{R} ^{2} \varphi ^{0},  } ^2\Bigr] ^{1/2} 
\end{equation*}

\begin{theorem}\label{t.H^1equiv}
All of the norms below are equivalent, and can be used as 
a definition of real $ \operatorname {Re}H ^{1} (\mathbb C_+^d)$. 
\begin{gather*}
\norm \operatorname M\cdots \operatorname M f.1.\,,
\qquad
\norm \operatorname S\cdots \operatorname S f.1.\,,
\qquad
\sum _{\sigma \in \{0,1\}^d} \norm \operatorname P _{\sigma } f.1.\,,
\\
\sum _{j=1} ^{d} \sum _{\operatorname A_j\in \{\operatorname I, \operatorname H_j\} } 
\Norm \prod _{j=1} ^{d} \operatorname A_j f.1.\,.
\end{gather*}
In the last line, we are summing over all choices of operators $\operatorname  A_j$ being 
either the identity operator, or $ \operatorname H_j$, the Hilbert transform 
computed in the $ j$th direction. 
\end{theorem}

\subsection{$ \textup{BMO} (\mathbb C_+^d)$}

The dual of the real Hardy space is $\operatorname {Re}\operatorname H^1(\mathbb C_+^d)^\ast=\text{BMO}(\mathbb C_+^d)$, the $d$--fold product $\text{BMO}$ space. It is a Theorem of S.-Y.~Chang and R.~Fefferman 
\cite{MR82a:32009} that this space has a characterization in terms of the product Carleson measure introduced above. 
We need the product wavelet basis.   For a rectangle $R=\prod_{j=1}^d R_{(j)}\in\mathcal D^d$
set 
 \begin{equation*}
w_R(x_1,\ldots,x_d)=\prod_{j=1}^d w_{R_{(j)}}(x_j)=\operatorname {Dil}_R ^{2} w _{[0,1]^d}(x)
\end{equation*}
The basis $\{ w_R\mid R\in\mathcal D^d\}$ is the $d$--fold tensor product of the wavelet
basis. 
We use the same notation $ w_R $ and $ v_R$ for the Meyer wavelet basis, and the analytic 
Meyer wavelet basis.  
Define 
\begin{equation} \label{e.BMOdef}
\lVert b\rVert_{\text{BMO} (\mathbb R^{d})}\simeq{}
\sup_{ U\subset\mathbb R^d} \Bigl[\abs{   U}^{-1} 
\sum_{R\subset U}\abs{\ip b,w_R,}^2\Bigr] ^{1/2} 
\end{equation}
where we have replaced the Haar wavelets  by the Meyer wavelets on the right.

It is the Theorem of Chang and Fefferman that 
\begin{theorem}\label{t.changfefferman}
We have the equivalence of norms 
\begin{equation*}
\norm f. (\operatorname {Re}H^1 (\mathbb C_+ ^{d})) ^{\ast}. 
\simeq \norm f. \textup{BMO}(\mathbb R^d). 
\end{equation*}
That is, $ \textup{BMO}(\mathbb R^d)$ is the dual to $ \operatorname {Re}H^1 (\mathbb C_+ ^{d})$.
\end{theorem}

To define \emph{analytic}
$ \textup{BMO} (\mathbb C_+^d)$, it suffices to replace the Meyer wavelets above by the analytic 
Meyer wavelets.

\subsection{Journ\'e's Lemma}

The explicit definition of $ \textup{BMO}$ in (\ref{e.BMOdef}) is quite difficult 
to work with. In the first place, it is not an intrinsic definition, in that 
one needs some notion of wavelet to define it.  Secondly, the supremum is 
over a very broad class of objects: All subsets of $ \mathbb R ^{d}$ of finite measure. 
There are simpler definitions,
(that unfortunately are not intrinsic) that in particular circumstances are 
sufficient.  

For our purposes, there are two appropriate definitions.  
Set $ \norm f.\text{BMO}(\text{rec}).$ to be the supremum in (\ref{e.BMOdef}), but 
with the important restriction that the sets $ U$ are taken to be rectangles.  
Historically, this was the first natural guess for the correct definition of 
$ \textup{BMO} (\mathbb C_+^d)$.  But, in a key moment, L.~Carleson 
\cite{carleson-example} produced examples of 
 functions which acted as linear functionals on $\operatorname H^1 (\mathbb C_+^d)$ with norm one, yet 
had arbitrarily small $\operatorname{BMO}(\text{rec})$ norm.  This example is recounted at the beginning of 
R.~Fefferman's article \cite{MR81c:32016}. 

Despite this fact, Journ\'e Lemma shows that in certain circumstances, the 
rectangular $ \text{BMO}(\text{rec})$ norm can dominate the $ \textup{BMO} (\mathbb C_+^d)$ 
norm.  
Let us state this Lemma in the case of $ \mathbb C_+^2$ before moving to the 
more sophisticated variants that we will need in three and higher parameters.

Given a set $ U\subset \mathbb R ^2$ of finite measure, let 
\begin{gather*}
\operatorname {Emb}(R; U)=\sup \{ \mu >1\mid 
\mu R\subset V\,\}
\\
V \eqdef \{ \operatorname M \operatorname M  \mathbf 1 _{ \{ \operatorname M \operatorname M \mathbf 1 _{U}>\tfrac12\} }
>\tfrac12\}
\end{gather*}
This is defined for rectangles $ R\subset U$, where $ \mu R$ denotes the rectangle 
with the same center as $ R$, which is dilated by a factor of $ \mu $ in \emph{ all} directions. 
Notice that we have $ \abs{ V} \lesssim \abs{ U}$, and that $ V$ is a natural `dilate of 
$ U$.'  The function $ \operatorname {Emb}(R; U) $ is a measure of how deeply embedded 
$ R$ is inside of $ U$.

A key distinction in two and higher parameters concerns 
collection of rectangles $ \{R\mid R\subset U,\ \operatorname {Emb}(R; U)\simeq \mu \}$. 
This collection of rectangles is \emph {not pairwise disjoint}, but their 
overlap is, in appropriate sense, at worst logarithmic in $ \mu $.  
A  formulation of this principle is easiest in two parameters.

\begin{lemma}[Journ\'e's Lemma \cite {MR87g:42028} in $ \mathbb R ^{2}$]\label{l.journe2} 
For any function $ f$, and $ \epsilon >0$, we have the inequality below valid 
for all sets $ U\subset \mathbb R ^2$ of finite measure.
\begin{equation*}
\NOrm \sum_{R\subset U} \operatorname {Emb}(R;U) ^{-\epsilon }
\ip f, w_R, w_R . \textup{BMO} (\mathbb R^d). 
\lesssim \norm f. \textup{BMO}( \textup{rec}).
\end{equation*}
The implied constant depends only on $ \epsilon >0$. 
\end{lemma}

Notice that the last inequality is that the $ \textup{BMO} $ norm is 
dominated by the (generally smaller) $ \textup{BMO(rec)}$ norm.
Carleson's examples show that this inequality is false if we do not `dampen' the 
wavelet coefficients in some way.  Journ\'e's insight is that this can be done 
with the geometric notion of the enlargement term.

We will need this observation in the case of $ \mathbb C_+^2$.  But, the rectangular 
norm is ill suited to our needs in three and higher dimensions.  We make this 
definition, which reduces to rectangular $ \textup{BMO}$ in dimension $ 2$. 

Say that a collection of rectangles $ \mathcal U\subset \mathcal D ^{d}$
\emph{ has $ d-1$ parameters} iff there is a choice of coordinate $ j$ so that 
for all $ R,R'\in \mathcal U$ we have $ R _{(j)}=R _{(j)}'$, that is the $ j$th coordinate 
of the rectangles agree.

We then define 
\begin{equation}\label{e.BMOd-1} 
\norm f. \textup{BMO} _{-1} (\mathbb R^d). 
\eqdef \sup _{\substack{\textup{$ \mathcal U$ has $ d-1$ } \\ \textup{parameters} }} 
\Bigl[\abs{ \operatorname {sh} (\mathcal U)} ^{-1} \sum_{R\in \mathcal U} 
\abs{ \ip f,w_R,} ^2 \Bigr] ^{1/2} 
\end{equation}
A collection of rectangles has a \emph {shadow} given by 
$ \operatorname {sh} (\mathcal U) \eqdef \bigcup\{R\mid R\in \mathcal U\} $. 
Observe that in $ d=2$ this reduces to
the rectangular $\textup{BMO}$ definition.  We use the $ -1$ subscript to indicate that 
we have `lost one parameter' in the definition.

The extension of Journ\'e's Lemma that we need replaces the $ \textup{BMO} (\textup{rec})$ 
norm by this $ \textup{BMO} _{-1} (\mathbb C_+^d)$ norm.  Yet one more refinement 
is essential for our needs, that the `dilate' of the set $ U$ be taken with considerably more care, 
and in particular should be just a little bit bigger than $ U$ in measure, see (\ref{e.d-1B}).

\begin{lemma}[Journ\'e's Lemma in $ d-1$ parameters]\label{l.journed-1}
For all $ \eta >0$, and  collections of rectangles  $ \mathcal U$
whose shadow has finite  measure, we can 
construct $ V\supset \operatorname {sh}(\mathcal U)$ and a function 
$ \operatorname {Emb}\mid \mathcal U\longrightarrow[1,\infty )$ so that 
\begin{gather} \label{e.d-1A}
\operatorname {Emb} (R)\cdot R\subset V,\qquad R\in \mathcal U\,,
\\ \label{e.d-1B}
\abs{ V}<(1+\eta ) \abs{ \operatorname {sh} (U)}\,,
\\ \label{e.d-1C}
\NOrm \sum_{R\subset U} \operatorname {Emb}(R;U) ^{-2d }
\ip f, w_R, w_R . \textup{BMO} (\mathbb R^d). \le
K_\eta  \norm f. \textup{BMO} _{-1}(\mathbb R ^{d}).
\end{gather}
The last inequality holds for all functions $ f$, with the constant $ K _{\eta }$ 
depending only on $ \eta $. 
\end{lemma}

Notice that the power on the embeddedness term   in (\ref{e.d-1C})
is quite large, twice the number of parameters.
Also, concerning the conclusions, if we were to take $\operatorname {Emb}(R)\equiv1$, then 
certainly the first conclusion (\ref{e.d-1A}) would be true.  But, the last 
conclusion would be false for the Carleson examples in particular.  This choice is 
obviously not permitted in general. 

The formulations of Journ\'e's Lemma given here are not the typical ones 
found in Journ\'e's original Lemma, or  J.~Pipher's extension to three dimensional 
case.  These papers give the more geometric formulation of these Lemmas, 
and J.~Pipher's article implicitly contains the geometric formulation 
needed to prove the Lemma above (provided one is satisfied with 
the estimate $ \abs{ V} \lesssim \abs{ \operatorname {sh} (\mathcal U)}$).  See 
Pipher \cite{MR88a:42019}. 
Lemma~\ref{l.journed-1}, as formulated above, was found in Lacey and Terwelleger \cite{witherin}; 
the two dimensional variant (which is much easier) appeared in Lacey and Ferguson 
\cite{MR1961195}.  The paper of Cabrelli, Lacey, Molter and Pipher 
\cite{math.CA/0412174} is a comprehensive survey of issues related to Journ\'e's Lemma. 
See in particular Sections 2 and 4.
We refer the reader to it for more information on this subject.

   \section{Multiparameter Paraproducts}  \label{s.multi}
 
We now consider paraproducts formed over sums of dyadic rectangles in  $\mathbb R^d$.
Let us say that a function $\varphi $  \emph { is adapted to a rectangle $R= \mathop\otimes _{j=1}^dR _{j}$ } iff 
$ \varphi (x_1,\dots ,x_d)=\prod _{j=1}^d\varphi_j(x_j) $, with each  $\varphi_j $ adapted to the interval $R_j $ 
in the sense of (\ref{e.adapted}).

Our paraproducts are of the same general form
\begin{equation*}
\operatorname B(f _{1},f _{2})\eqdef{}\sum _{R\in\mathcal R}
\frac{\varphi _{3,R}}{\abs{R}  ^{\frac{1}2}}\prod _{v=1}^{2} 
\ip f _{v},\varphi _{v,R},.
 \end{equation*}
 Here, we let $\mathcal R \eqdef{}\mathcal D^d $ be the class of dyadic rectangles. 
 
 The Theorem in this setting is 
 
 \begin{theorem}\label{t.2}  Let  $1<p < \infty $, and  $ J\subset \{1,\dotsc,d\}$.
 Assume that for each choice of coordinate $1\le j \le d$,  and $ v=1$ 
 \begin{equation}  \label{e.2-zeros}
 \int _{\mathbb R }\varphi _{v,R }(x_1,x_2,\dotsc,x_n)\; dx_j=0,\qquad \text{for all $x_ k$ with $k\not=j$ and all $R $ }.  
 \end{equation}
 In addition, for each $ 1\le j\le d$ and all $ R\in \mathcal D ^{d}$,
 assume that the condition above holds for $ \varphi _{v,R}$, where 
 $ v=2$ if $ j\in J$ and $ v=3$ if $ j\not\in J$.  
 Then, we have the inequality 
 \begin{equation}  \label{e.B-2-}
 \operatorname B\mid  \textup{BMO}(\mathbb C_+^d)\times L ^{p}\longrightarrow L^ {p} .
 \end{equation}
 \end{theorem}  

 We are not stating this result in greatest generality.  It was first discussed in the 
 the paper of Journ\'e \cite {MR88d:42028}.  Recently, the result has received 
 new attention, and extension.  See Muscalu, Pipher, Tao and Thiele \cites{camil,math.CA/0411607}.
 Our discussion is drawn from Lacey and Metcalfe \cite{math.CA/0502334}. 
And in particular, this last paper proves this Theorem.

  The critical distinction comes from the assumption about the zeros, (\ref{e.2-zeros}).  
  We need zeros in every coordinate  on the functions that land on 
 the  $\textup{BMO} (\mathbb C _+ ^{d})$ function.  There is one more zero 
 in each coordinate, and they can be split up between the second and third functions.

 Notice that there are many different types of paraproducts.  The first case, with the greatest similarity to the one parameter case, is where 
 we have, for example,  $x_j $ zeros in first and second positions for all $1\le j\le d$.   
 The other cases do not have a proper analog in the one parameter case.

 \bigskip 
 
 We will have need of paraproducts which are presented in a somewhat different 
 way, in analogy to Theorem~\ref{t.U}.  We make some definitions.  For 
 $\vec s\in \mathbb Z ^{d}$, let us set 
 \begin{equation*}
\Delta \operatorname U _{\vec j} =\sum_{\substack{R\in \mathcal D ^{d}\\ 
\abs{ R_s}=2 ^{j_s}\,,\ 1\le s \le d}} u _{R} \otimes u_R\,. 
\end{equation*}
 
For a subset of coordinates  $ J\subset \{1,\dotsc,d\}$ set 
\begin{equation*}
\operatorname U _{\vec \jmath,J} \eqdef \sum_{\substack{\vec k\in \mathbb Z ^{d}\\
k_s=j_s\,, \ s\in J \\ k_s\ge j_s\,, s\not\in J}}
\Delta \operatorname U _{\vec k}
\end{equation*}
For those coordinates $ s\in J$, we take the wavelet projection onto 
that scale, while for those coordinates $ s\not\in J$, we sum over all 
larger scales. 

Write $ R' \lesssim _{J} R $ iff  $ \abs{ R' _{s}}\le \abs{ R _{s}}$ for 
$ s\not\in J$ and $ \abs{ R'_s}=\abs{ R_s}$ for $ s\in J$. 

\begin{theorem}\label{t.Ud} For all $ J\subset \{1,\dotsc,d\}$, 
and $ \vec k\in \mathbb Z ^{d}$ with $ \norm \vec k.\infty .\le 8$, 
we have 
\begin{equation} \label{e.Ud}
\Norm \sum _{\vec \jmath\in \mathbb Z ^{d}} (\Delta 
\operatorname U _{\vec \jmath, J}\, b)\cdot 
\overline{\operatorname U _{\vec \jmath+\vec k,J} \varphi } .2. 
\lesssim \norm b.\textup{BMO}(\mathbb C_+ ^{d}). \norm \varphi .2.
\end{equation}
Moreover, suppose we have the following separation condition: 
Fix an integer $ A>0$. 
Suppose that 
\begin{equation}\label{e.Localized}
\textup{if $ \ip b,u_ {R'},\neq0$, $ \ip \varphi , u _{R},\neq0$  with 
$ R' \lesssim _{J} R$, then $ AR \cap R'=\emptyset $. }
\end{equation}
We then have the estimate 
\begin{equation} \label{e.localized}
\Norm \sum _{\vec \jmath\in \mathbb Z ^{d}} (\Delta  \operatorname 
U _{\vec \jmath, J} b)\cdot 
\overline{\operatorname U _{\vec \jmath+\vec k,J} \varphi } .2. 
\lesssim A ^{-100d}\norm b.\textup{BMO}(\mathbb C_+ ^{d}). \norm \varphi .2.
\end{equation}
Implied constants are independent of the choice of $ \vec k$. 
\end{theorem}

\begin{proof}
The method of proof is quite similar to that of Theorem~\ref{t.U}.  

We treat a special case with a brute force approach.  Consider 
\begin{equation*}
\sum_{_{\vec \jmath\in \mathbb Z ^{d}}}(\Delta \operatorname U _{\vec \jmath} b )
\cdot 
(\overline{\Delta \operatorname U _{\vec \jmath+\vec k} \varphi })\,, \qquad k\in \mathbb Z ^{d},\ 
\norm \vec k.\infty .\le8\,.
\end{equation*}
We claim that the two estimates of the Theorem hold for these operators.  

Fix an integer $B \ge2$
Let $ \pi \mid \mathcal D ^{d}\longrightarrow \mathcal D ^{d}$ be a map so that 
for all $ R\in \mathcal D ^{d}$ we have $ 2 ^{k_s}\abs{ R_s}=\abs{ \pi (R)_s}$
for $ 1\le s\le d$.  In addition, the distance between $ R$ and $ \pi (R)$ is 
essentially constant.  Namely, 
\begin{equation*}
B ^{-1}\le \operatorname M\cdots \operatorname M \mathbf 1 _{R}(c (\pi (R))\le (B -1)^{-1}
\end{equation*}
(In the current setting, it is most natural to use the maximal function to 
measure distances.) 
Then, the sum 
\begin{equation*}
\sum _{R\in \mathcal D ^{d}} 
\frac {\ip b, u_R,} {\sqrt \abs{ R}} \overline{\ip \varphi ,u_ {\pi (R)}, }  \cdot 
\sqrt {\abs{ R}}\, u_R \, \overline{u _{\pi (R)}}
\end{equation*}
is a paraproduct, with zeros in all coordinates for both $ b$ and $ \varphi $.  
This is not necessarily the case of the third place, but $B ^{200d}\sqrt {\abs{ R}}\, u_R \, \overline{u _{\pi (R)}} $
is adapted to $ R$ with constant independent of $ B$.
To prove (\ref{e.Ud}), we then sum over $ B\ge1$; to prove (\ref{e.Localized}), 
sum over $ B\ge A$.

\bigskip 

We recall the `father wavelet'  $ W$ from the proof of Theorem~\ref{t.U}.  
For a subset of coordinates $ J\subset \{1,\dotsc,d\}$ we set 
\begin{equation*}
W _{R,J}(x_1,\dotsc ,x_d) \eqdef \prod _{s\in J} u _{R_s}(x_s)\cdot 
\prod _{s\not\in J} W _{R_s}(x_s)
\end{equation*}
thus, in the coordinates in $ J$ we take an analytic Meyer wavelet, and 
for those coordinates not in $ J$ we take a father wavelet. 

Observe that 
\begin{equation*}
\operatorname U _{\vec \jmath, J} =\Bigl[ \sum _{\substack{R\in \mathcal D ^{d}\\
\abs{ R_s}=2 ^{j_s}}} W _{R,J} \otimes W _{R,J}\Bigr] \circ \operatorname P _{\oplus}.
\end{equation*}
For $ \vec 9=(9,\dotsc,9)$, we need to provide the two bounds for the Theorem for 
\begin{equation*}
\sum _{\vec \jmath\in \mathbb Z ^{d}} 
(\Delta \operatorname U _{\vec \jmath} b)\cdot \overline{\operatorname U _{\vec \jmath +\vec 9,J} \operatorname P _{\oplus}}
\end{equation*}

For an integer $ A$, and map $ \pi $ as above, it suffices to consider the sum 
\begin{gather*}
\sum_{R\in \mathcal D ^{d} } 
\frac {\ip b, u_R,} {\sqrt {\abs{ R}}} \, \overline{\ip \varphi , W _{ \pi (R),J} ,}  
\, \psi _{R,J}\,,
\\
\psi _{R,J} \eqdef \sqrt {\abs{ R}}\, u_R \cdot \operatorname W _{\pi (R),J}\,.
\end{gather*}
This is a paraproduct.  Note that for the function $ b$, we have zeros in all 
coordinates; for the function $ \varphi $, we have zeros in coordinates $ s\not \in J$; 
the function $ \psi _{R,J}$ has zeros in those coordinates $ s\in J$. Finally, 
$
B ^{200d} \psi _{R,J} 
$ 
is adapted to $ R$ with constants that are independent of $ B$ or the choice of $ \pi $. 
Thus, the two claimed inequalities of the Theorem hold for these sums,  and this completes 
the proof of the Theorem. 

\end{proof}

\section{Nehari Theorem in Several Variables} 

The Hankel operators we are are concerned with are maps from 
$ H^2 _{\oplus } (\mathbb C_+^{d})$ to $ H ^2 _{\oplus} (\mathbb C_+^{d})$ given by 
\begin{equation}\label{e.dHankel}
\operatorname H _{b} \varphi  \eqdef \operatorname P _{\oplus}  \operatorname M _{b} 
\overline \varphi \,. 
\end{equation}
This definition only depends upon $ \operatorname P _{\oplus } b$.

\begin{remark}\label{r.little} These are the `little' Hankel operators, in that 
we are taking the `smallest' reasonable projection above.  To define the `big' 
Hankel operators, one would replace $\operatorname P _{\ominus} $
 above by $\operatorname I- \operatorname P _{\oplus }$.  
 We refer the reader to Cotlar and Sadosky \cite {MR1284610} for the 
 theory of these `big' Hankel operators. 
\end{remark}

The Nehari Theorem in this context is: 

\begin{multinehari}[\cites{MR1961195,witherin}]We have the equivalence 
\begin{equation}\label{e.neharid}
\norm \operatorname H _{b }.. \simeq \norm \operatorname P _{\oplus} b. \textup{BMO} (\mathbb C_+^{d}).
\end{equation}
where the latter space is S.-Y.~Chang and R.~Fefferman $ \textup{BMO}$, the dual 
to the Hardy space $ H ^1 (\mathbb C_+ ^{d})$.
\end{multinehari}

This theorem has equivalent statements; the most obvious of these concerns the 
multiparameter commutator 
\begin{equation*}
\operatorname C (b,f) \eqdef [\cdots[ \operatorname M_b, \operatorname H_1],\dotsc, \operatorname H_d]
\end{equation*}
where $ \operatorname H _{j}$ denotes the Hilbert transform computed in the $ j$th coordinate. 

Less obviously, there is an an equivalent formulation in terms of factorization.  As we have 
commented, the classical factorization of $ H^1$ functions given in Proposition~\ref{p.factor} 
does not extend to $ H^1 (\mathbb C_+^{d})$.  The Nehari theorem is equivalent to 
\emph {weak factorization}.  The formalization of this is done in terms of 
a tensor products of $ H^2 (\mathbb C_+^{d})$.

We define a projective tensor product norm by 
\begin{equation*}
\norm f.H^2 (\mathbb C_+^{d})\widehat \otimes H^2 (\mathbb C_+^{d}). \eqdef 
\inf \Bigl\{  \sum_{j} \norm \varphi _j.H^2 (\mathbb C_+^{d}).\norm \psi .H^2 (\mathbb C_+^{d}).
\mid f=\sum _{j} \varphi _j \psi _j \,,\ \varphi _j, \psi _j \in H^2 (\mathbb C_+^{d}) \Bigr\}
\end{equation*}

\begin{theorem}\label{t.equiv} Any one of   equivalences of norms below are 
consequences of the other equivalences. 
\begin{align}\label{e.nehariNorm}
\norm \operatorname H_b..& \simeq   \norm \operatorname P _{\oplus }
b. \textup{BMO} (\mathbb R^{d}).
\\  \label{e.commd}
\norm \operatorname C_b.2\to 2.& \simeq   \norm   
b. \textup{BMO} (\mathbb C_+^{d}).
\\ \label{e.weakd}
\norm f.H^1 (\mathbb C_+^{d}).
& \simeq 
\norm f.H^2 (\mathbb C_+^{d})\widehat \otimes H^2 (\mathbb C_+^{d}). 
\end{align}
In the first two, we take the $ \textup{BMO} (\mathbb R ^{d})$ norm to be real valued $ \textup{BMO}$. 
In the second two, $ \textup{BMO}(\mathbb C _+ ^{d})$ is analytic $ \textup{BMO}$. 
\end{theorem}

The last equivalence of norms is the weak factorization statement in $ H^1 (\mathbb C_+^{d})$. 
It explains in part why the factorization proof of the one parameter Nehari theorem is 
so easy: The factorization property is stronger than Nehari's Theorem.

\begin{corollary}\label{c.nehari} We have 
\begin{equation*}
\norm \operatorname H _{b}..=\inf \bigl\{ \norm \beta .\infty . 
\mid \operatorname P_\oplus b=\operatorname P_\oplus \beta \bigr\}\,.
\end{equation*}
\end{corollary}

Theorem~\ref{t.equiv} was known and elementary;  once weak factorization 
(\ref{e.weakd}) is known, Corollary~\ref{c.nehari} is easy.  Thus,
 the Multiparameter Nehari Theorem is the main point. x
The   inequality 
 $ \norm \operatorname H _{b}.. \lesssim \norm b. \textup{BMO} (\mathbb C_+^{d}).$, 
turns out to be quite easy---it is a consequence of the trivial 
inclusion in the weak factorization statement. 
The issue is to establish the lower bound on the norm of the Hankel operator. 

The central difficulty here lies in the subtle nature of $ \textup{BMO}$ in 
the higher parameter case. 
 The proof we give is 
 an induction on $ d$, using weak factorization in $ H^1 (\mathbb C_+^{d-1})$ 
 in a critical moment.    Appealing to weak factorization will give us 
 a lower bound in terms of $ \textup{BMO} _{-1}$.  And so we need to `bootstrap' 
 from this weaker inequality to the stronger inequality.  The boostrapping argument appeals 
 to the Journ\'e Lemma.

It suffices to assume that $ b=\operatorname P _{\oplus } b \in \textup{BMO}(\mathbb C_+^{d})$ 
is of norm one, and find an absolute lower bound on $ \norm H _{b}..$. 
We begin by using the induction hypothesis to establish 
\begin{equation*}
\norm H _{b}.. \gtrsim \norm b. \textup{BMO} _{-1} (\mathbb C_+^{d}).,
\end{equation*}
  where the latter norm is   $ \textup{BMO}$ norm  `with one less parameter' 
defined in (\ref{e.BMOd-1}).
 Thus, we are free to impose the additional hypothesis that 
 $ \norm b. \textup{BMO} _{-1} (\mathbb C_+^{d}).$ is less than some 
 fixed, absolute constant. 
Observe that implicitly, this forces $ b$ to be the type of functions which 
Carleson discovered.

  Yet, Journ\'e's Lemma gives modest sufficient conditions for this 
  impoverished norm to dominate the true $ \textup{BMO}$ norm.  The lower bound 
  for the norm of $ \operatorname H _b$ can then be explicitly estimated 
  as a main term, plus several error terms.  Each of the error terms is 
  a paraproduct, which can be controlled with Journ\'e's Lemma and the 
  fact that the improvised norm is small.

\begin{proof}[Proof of Theorem~\ref{t.equiv}.]
We discuss the proof of Theorem~\ref{t.equiv} and Corollary~\ref{c.nehari}. 
Observe that the computation (\ref{e.dualNorm}) is quite general. In the 
language we have introduced above, it shows immediately that  
\begin{equation} \label{e.nehari*}
\norm \operatorname H _{b}.. 
\simeq \norm \operatorname P _{\oplus } b. 
H^2 (\mathbb C_+^{d})\widehat \otimes H^2 (\mathbb C_+^{d}) ^\ast. 
\end{equation}
That is, the Hankel norms are equivalent to the \emph {dual norm} of the 
tensor product norm. 

The equivalence of (\ref{e.nehariNorm}) and (\ref{e.weakd}) is then immediate. 

Concerning the commutator, and (\ref{e.commd}), as in the one parameter case, 
the commutator is seen to be a sum of $ 2^d$ Hankel operators.  
Indeed, for $ \sigma \in \{-,+\}^d$, consider the composition 
$
\operatorname C_b  \operatorname P _{\sigma }
$
In the definition of the commutator, we are free to replace 
the $ j$th Hilbert transform $ \operatorname H_j$ by $ \operatorname P _{-\sigma (j),j}$, 
since $ H_j=\pm(\operatorname I-2  \operatorname P _{-\sigma (j),j})$, and 
the identity commutes with everything.  Thus, 
\begin{equation*}
\operatorname C_b \operatorname P _{\sigma }=\pm 2^d 
\operatorname P_ {-\sigma } \operatorname M_b \operatorname P _{\sigma }\,. 
\end{equation*}
From this, the equivalence of (\ref{e.nehariNorm}) and (\ref{e.commd}) is immediate. 

\end{proof}

\begin{proof}[Proof of Corollary~\ref{c.nehari}.]
We can assume that the symbol of the Hankel operator $H _{b} $ is in analytic $ \textup{BMO}$. 
Then, (\ref{e.nehari*}) and (\ref{e.weakd}) show that $ b$ defines a bounded 
linear functional on $ H^1 (\mathbb C_+^{d})\subset L^1(\mathbb R ^{d})$.  
Appeal to the Hahn Banach Theorem to extend this linear functional to all of $ L^1 (\mathbb R ^{d})$, 
with the same norm.  The Corollary follows. 
\end{proof}

\section{Proof of Multiparameter Nehari Theorem} 

The upper bound on the norm of a Hankel operator is easy.  Observe that, trivially, 
\begin{equation*}
 H^2 (\mathbb C_+^{d})\widehat \otimes H^2 (\mathbb C_+^{d})
 \subset H^1 (\mathbb C_+^{d})\,.
\end{equation*}
For the dual spaces, we have the reverse inclusion. 
In particular, the $ \textup{BMO}(\mathbb C_+^{d})$ norm is larger than 
the dual tensor product norm. 
Thus, by (\ref{e.neharid}), 
\begin{align*}
\norm \operatorname H_b.. 
&\simeq \norm \operatorname P _{\oplus } b. 
H^2 (\mathbb C_+^{d})\widehat \otimes H^2 (\mathbb C_+^{d}) ^\ast. 
\\
& \lesssim \norm \operatorname P _{\oplus } b. \textup{BMO}(\mathbb C_+^{d}).
\end{align*}
Thus, the primary difficulty is in establishing the lower bound on the norm 
of the Hankel operator.

\subsection{The Initial Lower Bound}

The proof is by induction on dimension $ d$, and we take the classical Nehari 
Theorem as the base case in the induction.  Thus, we assume that (\ref{e.neharid}) 
holds in dimension $ d-1\ge1$, and prove it in dimension $ d$. 

Take $ b$ to be in  analytic $ \textup{BMO}(\mathbb C_+^{d})$, 
and of norm one.  We recall that this means in particular, that we have 
\begin{equation}\label{e.bmoForb}
\sup _{U\subset \mathbb C_+^{d}} \abs{ U} ^{-1} 
\sum _{\substack{R\in \mathcal D ^{d}\\ R\subset U }} 
\abs{ \ip b,v_R,}^2=1\,,
\end{equation}
where we recall that the supremum is over all subsets $ U$ of finite measure, and 
that the functions $ v_R$ are the analytic Meyer wavelets associated to dyadic 
rectangles in $ \mathbb C_+^{d}$.

Let us argue that
\begin{equation}\label{e.bmod-1}
\norm \operatorname H_b .. \gtrsim \norm b. \textup{BMO}_{d-1} (\mathbb C_+^{d}).
\end{equation}
This last norm is given in (\ref{e.BMOd-1}), and in particular, it is a supremum as in 
(\ref{e.bmoForb}), with an additional restriction on rectangles that 
contribute to that sum.  

Now, this inequality we are to prove, by (\ref{e.nehari*}), reduces to showing 
\begin{equation}\label{e.bmodd}
\norm \operatorname P _{\oplus } b.
(H^2 (\mathbb C_+^{d})\widehat \otimes H^2 (\mathbb C_+^{d}))^\ast. 
\gtrsim \norm \operatorname P _{\oplus } b. \textup{BMO}_{d-1} (\mathbb C_+^{d}).\,. 
\end{equation}
We can assume that $ b=\operatorname P _{\oplus } b$ is a Schwartz function,  
and that $ \norm b. \textup{BMO} _{d-1} (\mathbb C_+^{d}).=1$.  Thus, after 
a permutation of coordinate and a possible dilation,
we can take a collection of rectangles $ \mathcal U$ which achieves the 
 supremum in the $\textup{BMO} _{d-1} (\mathbb C_+^{d})$ norm. 
 
In particular, we can assume that 
\begin{itemize}
\item $ \abs{ \operatorname {sh}(\mathcal U)}=1$ ; 
\item there is an interval $ I$ of length one so that for all $ R\in \mathcal U$ we have 
$ R_1=I$; 
\item  for $ \psi =\sum_{R\in \mathcal U} \ip b,v_R, v_R$ we have $ \ip b,\psi ,=1$. 
\end{itemize}
Then, it suffices to see that $ \norm \psi .H^2 (\mathbb C_+^{d})\widehat \otimes H^2 (\mathbb C_+^{d}).
\lesssim  1$.

Write $ x=(x_1,x_2, \dotsc , x_d)\in \mathbb R ^{d}$ 
as $ (x_1,x')$ with $ x'=(x_2,\dotsc,x_d)\in \mathbb R ^{d-1}$.  
Each rectangle $ R\in \mathcal U$ has the same first coordinate. So the first 
coordinate in the in product that defines the Meyer analytic wavelet 
$ v_R$ is independent of $ R$. Therefore, 
we can write 
$ \psi (x)=\psi _{1}(x_1) \psi '(x')$  where $ \psi _1(x_1)\in H^1 (\mathbb R )$ 
is of norm one.  It can written as $ \psi _1=\alpha \cdot \beta $ 
with $ \alpha _1$ and $\beta _1$ of $ H^2(\mathbb R ) $ norm one.

$ \psi '$ satisfies something similar.  Observe that 
\begin{equation*}
\norm \psi ' . H^1 (\mathbb C_+^{d-1}). \le \abs{ U} ^{1/2} \norm \psi '.2.\le1\,.
\end{equation*}
Hence, $ \psi '$ is in $ H^1 (\mathbb C_+^{d-1})$, and is of norm  at most one. 
In fact, it has norm comparable to one, since  by construction 
$
\norm \psi '. \textup{BMO}(\mathbb C_+^{d-1}).\le1
$ 
and $ \ip \psi ',\psi ',=1$.
Thus, by the 
induction hypothesis, we have 
\begin{equation*}
\norm \psi '. H^1 (\mathbb C_+^{d-1}).\simeq 
\norm \psi '.H^2 (\mathbb C_+^{d-1})\widehat \otimes H^2 (\mathbb C_+^{d-1}). \simeq 1\,.
\end{equation*}
Thus, $ \psi '$ can be written as a sum of products of $ \alpha' _j\cdot\beta' _j$ with 
\begin{equation*}
\sum_{j} 
\norm \alpha '_j. H^2 (\mathbb C_+^{d-1}). 
\norm \beta  '_j. H^2 (\mathbb C_+^{d-1}). \simeq 1\,.
\end{equation*}
But then, it is clear that we can write 
\begin{equation*}
\psi (x_1,x') =\sum_{j} \alpha  (x_1) \alpha '_j (x')\cdot 
\beta   (x_1) \beta  '_j (x')
\end{equation*}
and so $ \norm \psi . H^2 (\mathbb C_+^{d})\widehat \otimes H^2 (\mathbb C_+^{d}). \lesssim  1 $.

\subsection{The $ \textup{BMO}(\mathbb C_+^{d})$ lower bound.} 

Our task is to `bootstrap' from the weaker inequality (\ref{e.bmod-1}).  
Namely, for an absolute constant $ \eta _{-1}$ whose value is to be specified, 
it suffices to consider Hankel symbols $ b$ which satisfy 
$b= \operatorname P _{\oplus } b$; $ b$ is Schwartz function; 
$ \norm b. \textup{BMO} (\mathbb C_+^{d}).=1$; and 
$ \norm b. \textup{BMO} _{-1} (\mathbb C_+^{d}). < \eta _{-1}$.  
(The subscript $ {} _{-1}$ mimics our notation for the reduced parameter $ \textup{BMO}$ space.)

We show by direct computation that $ \norm \operatorname H_b.. \gtrsim 1$, 
namely we will apply the Hankel operator to a particular $ H^2 (\mathbb C_+^{d})$
function, and provide a lower bound on the norm of the image.  

Here is how we select the test function to apply the Hankel to.  Select a collection 
of rectangles $ \mathcal U$ which achieve the supremum in the definition of 
$ \textup{BMO}(\mathbb C_+^{d})$ norm.  Thus, 
\begin{equation*}
\sum_{R\in \mathcal U} \abs{ \ip b,v_R,}^2=\abs{ \operatorname {sh} (\mathcal U)}\,.
\end{equation*}
Moreover, we can, after taking an appropriate dilation, that $ 
\abs{ \operatorname {sh} (\mathcal U)}=1$, and that if $ R\subset \operatorname {sh} (
\mathcal U)$, then $ R\in \mathcal U$. 

The function we apply the Hankel to the wavelet projection of $ b$ onto the 
wavelets associated with $ \mathcal U$, 
$ \alpha =\sum_{R\in \mathcal U}\ip b,v_R, v_R$. 
Observe that 
\begin{align*}
\norm \operatorname H_b \alpha ..&=\norm \operatorname P _{\oplus } \abs{ \alpha }^2
.2.
\\
& \gtrsim \norm \abs{ \alpha }^2.2.
\\
&= \norm \alpha .4.^2
\\
&\simeq \NOrm 
\Bigl[\sum _{R\in \mathcal U} \frac {\abs{ \ip b,v_R,} ^2} {\abs{ R}} \mathbf 1 _{R} \Bigr]
^{1/2} .4.^2
\\
&\ge 
\Bigl[\sum _{R\in \mathcal U}  {\abs{ \ip b,v_R,} ^2}  \Bigr] ^{1/2} \simeq 1\,.
\end{align*}
Here, we are relying on the symmetry of the Fourier transform of positive functions; 
Littlewood Paley inequalities, to pass to the 
wavelet square function; that $ \mathcal U$ has shadow equal to one in measure, 
and that $ L^4$ norms dominate $ L^2$ norms on a probability space. 
Thus, we have $ \norm \operatorname  H _{\alpha } \alpha ..\ge \eta _0>0 $, for absolute $ \eta _0$.

This is in fact our main estimate.  Our task is to show that for $ \eta _{-1}$ sufficiently 
small, that we have 
\begin{equation}\label{e.2DO}
\norm \operatorname H _{b-\alpha } \alpha ..<\tfrac12 \eta _0\,.
\end{equation}
This can be done with the aid of Journ\'e's Lemma.  

Fix a second small parameter $ \eta _{\textup{J}}$ whose value will be specified below. 
(The subscript $  {} _{\textup{J}}$ is for `Journ\'e.')
Apply Lemma~\ref{l.journed-1}. 
There is a set $ V\supset \operatorname {sh} (\mathcal U)$ 
and a function $ \operatorname {Emb} \mid \mathcal U\longrightarrow [1,\infty )$ 
for which these conditions hold. 
\begin{itemize}
\item $ \abs{ V}<1+\eta _{\textup{J}}$; 
\item $  \operatorname {Emb} (R) R\subset V$ for all $ R\in \mathcal U$; 
\item $ \norm \widetilde \alpha .\textup{BMO}(\mathbb C_+^{d}). \le K _{\eta _{\textup{J}}} \eta _{-1}$
\end{itemize}
where in the last line, we have 
\begin{equation}\label{e.tildealpha}
\widetilde \alpha \eqdef \sum_{R\in \mathcal U}  \operatorname {Emb}(R) ^{-2d} 
\ip b,v_R, v_R\,. 
\end{equation}

We now decompose the symbol $ b$.  We have already defined $ \alpha $. Set 
\begin{equation}\label{e.beta}
\beta \eqdef \sum_{\substack{R\subset V\\ R\not\in \mathcal U } }\ip b,v_R, v_R\,. 
\end{equation}
Thus, these are the rectangles with are `close' to $ \mathcal U$, but 
not in it, as defined by 
the set $ V$.  Define $ \gamma $ by $ b=\alpha +\beta +\gamma $. 
To verify (\ref{e.2DO}), it suffices to show that 
\begin{align}\label{e.2doBeta}
\norm \operatorname H _{\beta  } \alpha ..& <K \eta _{\textup{J}} ^{1/4}\,,
\\ \label{e.2doGamma}
\norm \operatorname H _{\gamma   } \alpha ..& <K _{\eta _{\textup{J}}} \eta _{-1} \,. 
\end{align}
One then specifies $ \eta _{\textup{J}}$ so that the top line is no more than 
$ \tfrac14 \eta _0$. The constant $ K _{\eta _{\textup{J}}}$ that appears in the second 
line is absolute, so we can then fix $ \eta _{-1}$ sufficiently small to prove 
(\ref{e.2DO}). 

\medskip 

The inequality for $ \beta $ is easily available to us, by the particular form 
of the Journ\'e Lemma we are using. Observe first that 
\begin{align*}
1+\sum _{\substack{R\subset V\\ R\not\in \mathcal U }}
\abs{ \ip b,v_R,} ^2 =
\sum_{R\subset V}\abs{ \ip b,v_R,} ^2 
\le1+\eta _{\textup{J}}\,.
\end{align*}
Therefore, $ \norm \beta .2.\le{} \sqrt {\eta _{\textup{J}}}$.  On the other hand, 
the $ \textup{BMO}(\mathbb C_+^{d})$ norm of $ \beta $ is less than or equal to one. 
Thus, we have $ \norm \beta .4. \lesssim \eta ^{1/4}$. 
A Hankel operator is at worst a product, thus 
\begin{equation*}
\norm \operatorname H _{\beta  } \alpha .. \le 
\norm \beta .4. \norm \alpha .4.\le K \eta _{\textup{J}} ^{1/4}\,. 
\end{equation*}
So it remains to verify (\ref{e.2doGamma}).

\subsubsection*{An Initial Calculation}

We make an explicit computation of a Hankel operator, in a manner similar to 
(\ref{e.commCases}). Namely, restricting attention to one dimension, we have 
\begin{equation}\label{e.HIJ} 
\operatorname H _{v_I} \overline {v_J}=\operatorname P _{+} (v_I \overline {v_J})
= 
\begin{cases}
0& 8 \abs{ J}< \abs{ I} 
\\
\operatorname P_+ (v_I \overline {v_J}) 
& \abs{ I}\le 8 \abs{ J}\le 64 \abs{ I}
\\
v_I \overline {v_J} & \abs{ I}< 8 \abs{ J}\,.
\end{cases}
\end{equation}
This follows from the Fourier localization properties of the Meyer wavelet.  
The Fourier support of the product 
 $ v_I \overline {v_J}$ is given by the convolution of the Fourier supports, 
which are specified by (\ref{e.meyer}).  If $ J$ is much smaller than $ I$, the 
product is purely antianalytic, 
giving us the first case above.  In the third case, $  v_I \overline {v_J}$ is 
purely analytic. 

\bigskip 

We apply the observation above to the term $ \operatorname H _{\gamma } \overline{\alpha }$. 
This leads us to the conclusion that 
\begin{gather*}
\norm \operatorname H _{\gamma } \overline{\alpha }..=
\NOrm 
\sum _{(R,R')\in \mathcal A} \ip b, u_R,\, \overline{\ip \varphi ,u_R,} \, u_R\, \overline{u_R}\,,
.2.
\\
\mathcal A \eqdef \{ (R,R')\mid 
R\subset U\,,\ R'\not\subset V\,, \abs{ R_s'}\le 64 \abs{ R_s}\,, 1\le s\le d\}\,.
\end{gather*}
It is essential to observe that this last sum can be written as a finite sum of the 
paraproducts in Theorem~\ref{t.Ud}, applied to the functions $ \alpha $ and $ \gamma $. 
This sum varies of choices of $ \vec k$ with $ \abs{ \vec k}\le 6$, and
arbitrary $ J\subset \{
1,\dotsc,d\}$. (The subset $ J$ consists of those coordinates  $ s$ for which 
$ \abs{ R_s} =2 ^{k_s} \abs{ R'_s}$.)

We  use  Theorem~\ref{t.Ud} to provide an estimate of the 
$ L^2$ norm of the sum above  an absolute constant times  $ \eta _{-1}$. 
In particular, we want to use the more technical estimate (\ref{e.localized})
 to achieve this end.

We will need to decompose the collection $ \mathcal A$ into appropriate
parts to which this estimate applies.  That is the purpose of this definition.
For an integer $ n\ge 1$, take 
\begin{equation*}
\alpha _{n} \eqdef \sum _{\substack{R\subset U\\ 
2 ^{n-1}\le 
\textup{Enl} (R; U)\le 2 ^{n}}} \ip b, u_R, u_R
\end{equation*}
We claim that 
\begin{equation}\label{e.claim1}
\norm \operatorname  H _{\gamma } \alpha _{n}.. \lesssim 2 ^{-n} \eta _{-1}\,. 
\end{equation}
It follows from Lemma~\ref{l.journed-1} that we have the estimate 
\begin{equation}\label{e.alpha<}
\norm \alpha _{n}. \textup{BMO}(\mathbb C _{+} ^{d}).
\lesssim 2 ^{2d n} \eta _{-1}\,,
\end{equation}
indeed, this is the point of this definition. From other parts of the expansion of the Hankel operator, we need to 
find some decay in $ n$.

Nevertheless, from this estimate and the  upper bound on Hankel operator norms, we have 
the estimate 
\begin{equation*}
\norm 
\operatorname  H _{\gamma } \alpha _{n}.. \lesssim 
\norm b. \textup{BMO}(\mathbb C _{+} ^{d}). \norm \alpha _{n}.2.
\lesssim 2 ^{2d n} \eta _{-1}.
\end{equation*}
We use this estimate for $ n<20$, say.

For $ n\ge20$, $ R\in \mathcal  U$ with $2 ^{n-1}\le 
\textup{Enl} (R; \mathcal U)\le 2 ^{n} $, and rectangle $ R'$ with $ (R,R')\in \mathcal A$, 
it follows that we must have $ 2 ^{n-9}R\cap R'=\emptyset$. 
That is, (\ref{e.Localized}) is satisfied with the value of $ A$ in that 
display being $ A\simeq 2 ^{n}$ for $ n\ge20$.  Thus, we conclude that 
\begin{equation*}
\norm 
\operatorname  H _{\gamma } \alpha _{n}.. \lesssim 
2 ^{-50n} \eta _{-1},\qquad n\ge 20\,. 
\end{equation*}
This completes our proof of (\ref{e.claim1}), and the proof of the
 lower bound on the norm of Hankel operators.

  
\end{document}